\documentclass{article}
\usepackage[ journal = preprint,   %% Replace XXX by one of the above journal codes.
lang=british,   %% Change to `american' if you use American English.
]{ems-journal}

\usepackage{bm}
\usepackage[algoruled,boxed,lined]{algorithm2e}
\newtheorem{remark}{Remark}[section]
\newtheorem{theorem}{Theorem}[section]
\newtheorem{assumption}{Assumption}[section]
\usepackage{mathtools}

\usepackage{lineno}
\newcommand{\Rd}{ {\mathbb{R}^d}}

\newcommand{\R}{ \mathbb{R}}
\newcommand{\E}{\mathcal{E}}

\newcommand{\ve}{\varepsilon}
\newcommand{\eps}{\epsilon}
\renewcommand{\d}{\textup{d}}
\newcommand{\m}{\textup{m}}
\newcommand{\supp}{\textup{supp}}
\newcommand{\law}{\textup{Law}}

\newcommand{\rank}{\textup{rank}}
\newcommand{\best}{\textup{best}}
\newcommand{\Parp}{\mathscr{P}}
\newcommand{\Par}{P}
\newcommand{\coll}{\mathscr{C}}

\newcommand{\argmin}{\textup{argmin}}

\renewcommand{\k}{{k}}
\renewcommand{\m}{\textup{m}}
\newcommand{\sve}{\sqrt{\ve}}

\newcommand{\OX}{\overline{x}}
\newcommand{\OY}{\overline{Y}}
\newcommand{\OV}{\overline{V}}
\newcommand{\OB}{\overline{B}}

\begin{document}

%%FOR SUBMISSION
%\mainmatter

%------
% Insert the title of your paper and (if necessary)
% a short title for the running head.
%------
\title{Kinetic models for optimization: a unified mathematical framework for metaheuristics}
%\titlemark{A review of kinetic models for stochastic particle methods in optimization}

%------
% Insert full names of the authors.
% Add further authors as follows:
%  \emsauthor{2}{}{}
%  \emsauthor{3}{}{}
% etc.
% Abbreviate first names for the running head.
%------
\emsauthor{1}{Giacomo Borghi}{G.~Borghi}
\emsauthor{2}{Michael Herty}{M.~Herty}
\emsauthor{3}{Lorenzo Pareschi}{L.~Pareschi}

%------
% Use \authormark if the list of authors is too
% long for the running head: \authormark{A.~Doe et al.}
%------

%------
% Add one \emsaffil and one \email for each author.
% NOTE: The address does NOT appear in the paper.
% It will probably be printed in an appendix.
%------
\emsaffil{1}{Maxwell Institute for Mathematical Sciences and Department of Mathematics, School of Mathematical and Computer Sciences (MACS), Heriot-Watt University, Edinburgh, UK \email{g.borghi@hw.ac.uk}}
\emsaffil{2}{Institut f\"ur Geometrie und Praktische Mathematik, RWTH Aachen University, Aachen, Germany \email{herty@igpm.rwth-aachen.de}}
\emsaffil{3}{Maxwell Institute for Mathematical Sciences and Department of Mathematics, School of Mathematical and Computer Sciences (MACS), Heriot-Watt University, Edinburgh, UK; Department of Mathematics and Computer Science, University of Ferrara, Italy \email{l.pareschi@hw.ac.uk}}

%------
% Add MSC 2020 codes according to www.ams.org/msc/msc2020.html.
% Secondary codes (in square brackets) are optional.
%------

%------
% Add a list of keywords.
%------

%------
% Optional: dedication
%------
%\chapterdedication{Dedicated to ...}

%------
% Insert your abstract.
%------
\begin{abstract}
Metaheuristic algorithms, widely used for solving complex non-convex and non-differentiable optimization problems, often lack a solid mathematical foundation. In this review, we explore how concepts and methods from kinetic theory can offer a potential unifying framework for a variety of metaheuristic optimization methods. By applying principles from collisional and non-collisional kinetic theory, we outline how particle-based algorithms like Simulated Annealing, Genetic Algorithms, Particle Swarm Optimization, and Ensemble Kalman Filter may be described through a common statistical perspective. This approach not only provides a path to deeper theoretical insights and connects different methods under suitable asymptotic scalings, but also enables the derivation of novel algorithms using alternative numerical solvers. While not exhaustive, our review highlights how kinetic models can enhance the mathematical understanding of existing optimization algorithms and inspire new computational strategies.
\end{abstract}

%%FOR SUBMISSION
%\makecontribtitle
%%FOR ARXIV
\maketitle

%------
% INSERT THE BODY OF THE PAPER HERE (except
% acknowledgments, funding info and bibliography)
%------

\section{Introduction}

In many practical applications, we are often required to solve optimization problems of the form
\begin{equation}
x^\star \in \underset{x\in \Omega}{\argmin}\, \E(x)
\label{eq:pb}
\end{equation}
where $\Omega$ is a given search domain, which could be finite or infinite. When dealing with complex problems—such as those where the objective function is highly non-convex, non-differentiable, or given as a black box—classical mathematical programming techniques often fail or become impractical.

To address these challenges, one may turn to heuristic optimization strategies that do not rely on gradient information and aim to find satisfactory solutions within a reasonable time frame.
These strategies, known as \textit{metaheuristics}, are designed to explore the search space and refine candidate solutions iteratively. The algorithms typically involve a combination of stochastic and deterministic updates that help balance the exploration and exploitation of the objective landscape. Some of the most well-known methods in this category include Genetic Algorithms (GA) \cite{holland1975}, Simulated Annealing (SA) \cite{kirk1983sa}, Particle Swarm Optimization (PSO) \cite{kennedy1995}, and Ant Colony Optimization (ACO) \cite{dorigo2006ant}. As the names suggest, each algorithm is inspired by a natural phenomena which also serve for a better understanding of the mechanisms.

The applications of these metaheuristic optimization techniques are vast and varied. They extend across fields such as engineering, where they can optimize design processes \cite{heur_eng}; finance, for portfolio optimization and risk management \cite{heur_finance}; logistics, to enhance supply chain operations \cite{heur_supply}; and machine learning, for tuning hyperparameters in complex models \cite{heur_ml}. Moreover, these methods are increasingly being utilized in the realms of healthcare for optimizing treatment plans \cite{heur_health}, environmental science for resource management \cite{heur_env}, and artificial intelligence for solving complex decision-making problems \cite{heur_ai}.

However, despite their popularity and success, metaheuristics lack a rigorous mathematical foundation. Traditional analytical tools from mathematical programming are of little use, as might be expected. A more fruitful approach is to compare metaheuristics to models from statistical physics, where complex systems are often described in probabilistic terms due to complex interactions and dependencies among particles. In optimization, these "particles" are candidate solutions, and the interplay between them can be studied using probabilistic models to better understand their dynamics and evolution.

While connections between statistical physics and optimization have long been explored, particularly for Simulated Annealing and discrete optimization algorithms \cite{magnus1996,doye1998,kirk2006stochastic,zecchina2005,holley1988}, a more systematic treatment of metaheuristics using continuous models based on partial differential equations (PDE) is relatively recent \cite{pinnau2017consensus,grassi2021pso,borghi2024ga,albi2023kinetic,pareschi2024sa}. Advances in the mathematical analysis of active particle systems \cite{carrillo2010particle,pareschi13,degond2013kinetic,diez2022review2} have provided refined tools for modeling the complex interactions between particles (or candidate solutions) in metaheuristics. One prominent example of this is the Consensus-Based Optimization (CBO) algorithm \cite{pinnau2017consensus}, which lends itself naturally to a mean-field description via nonlinear and nonlocal PDE models.

As recent studies on CBO and related methods have shown \cite{carrillo2018analytical,fornasier2024consensusbased,jin2020convergence,huang2021meanfield,hoffmann2024meanfield, Forna1, Forna2, Forna3,Borghi1, Borghi2, Borghi3, Wolfram, carrillo1, carrillo2, claudia1, claudia2, jin1, jin2, jin3, kalise}, there are several key advantages to modeling metaheuristics through kinetic equations:
\begin{itemize}
\item The statistical models can be more easily analyzed via classical PDE techniques to study algorithm convergence properties.
\item  They provide a higher-level, unified understanding of the underlying mechanisms, which can be used to improve existing methods.
\item By using different numerical techniques to discretize the PDE models, we can derive novel algorithms that may be more efficient.
\item Kinetic models offer a common framework to classify and connect different metaheuristic methods, revealing relationships through asymptotic scalings.
\end{itemize}

The aim of this review is to present a unified mathematical framework for several popular metaheuristics based on their kinetic descriptions. For each algorithm, we begin by outlining the heuristic and then derive the corresponding kinetic model. In the case of single-solution algorithms, we focus on the probability distribution of a single particle, whereas for multi-particle algorithms, we employ the mean-field approximation, assuming the system exhibits chaotic behavior in the many-particle limit \cite{sznitman1991chaos,diez2022review1}. We also collect known results on the convergence of the PDE models towards solutions to \eqref{eq:pb}. Furthermore, drawing on classical asymptotic limits in kinetic theory, such as the grazing collisions limit in the Boltzmann equation \cite{desvillettes1992asymptotics}, we illustrate how different metaheuristics can be connected through appropriate parameter scalings.

In the next sections, we will focus on continuous optimization problems where $\Omega  =\Rd$ in \eqref{eq:pb}, and consider the following algorithms: Simulated Annealing (Section \ref{sec:SA}), Genetic Algorithm (Section \ref{sec:ga}), Particle Swarm Optimization (Section \ref{sec:pso}) and the Ensemble Kalman Filter (Section \ref{sec:ekf}). These algorithms are among the most widely used metaheuristics, with many of them being readily available in popular optimization software libraries. They provide a diverse yet coherent sample of heuristic strategies, highlighting the potential of kinetic theory to unify and advance our understanding of metaheuristic optimization.

%%%% SA
\section{Simulated Annealing}\label{sec:SA}

\subsection{The heuristic strategy}

Introduced independently by \cite{kirk1983sa} and \cite{cerny1985sa}, Simulated Annealing (SA) is an extension to the Monte Carlo algorithm of Metropolis et al. \cite{metropolis1953}. The Metropolis--Hasting sampling method is based on the construction of a Markov chain which converges to the Boltzmann--Gibbs distribution 
\begin{equation}
f^T_\infty(x)   := \frac{e^{-\E(x)/T}}{Z_T} \,, 
\label{eq:gibbs}
\end{equation}
with $Z_T =\int e^{-\E/T} \d x$ normalizing constant. Taking the point of view of statistical physics, SA considers $\E(x)$ to be the energy associated with the state $x$, in our case $x\in \Rd$. The Boltzmann--Gibbs distribution $f_\infty^T$ captures the behaviour of a physical system in thermal equilibrium at finite temperature $T>0$. As the temperature decreases, $f_\infty^T$ concentrates more and more around low-energy configuration.  This behavior is related to the Laplace principle, which states that for any probability density $f$ 
\begin{equation}
\lim_{T \to 0} - T \log \left(\int e^{- \E(x)/T}f(x)\d x  \right) = \inf_{x\in \supp(f)} \E(x)\,,
\end{equation}
see \cite{dembo2010}. 
In SA, therefore, the Metropolis--Hasting Markov chain is modified such that the temperature is decreased during the computation. The process resembles the annealing strategy used to grow a single crystal, where the substance is first melted and, then, the temperature is carefully lowered, spending more and more time in the vicinity of the freezing point  \cite{kirk1983sa}.

In the SA algorithm, a single solution $x_\k$ at  temperature $T_\k>0$ is updated at every step $k$ according to the following rules:
\begin{itemize}
\item Generate a candidate solution $y_{k+1}$ via a random walk step
\begin{equation*}
y_{k+1} = x_\k + \sigma_\k \xi_\k
\end{equation*}
with, for instance, $\xi_\k$ randomly sampled from the standard normal distribution. The exploration coefficient $\sigma_k$ is time-dependent, and a classical choice is to set $\sigma_k = \sqrt{2T_\k}$.
\item If the candidate solution is better then the previous one, $\E(y_{k+1})  < \E(x_\k)$, it is accepted and so $x_{k+1} = y_{k+1}$. Otherwise, the candidate solution may still be accepted, but non-deterministically according to the rule
\begin{equation*}
x_{k+1} = 
\begin{cases}
y_{k+1} & \textup{with probability} \quad e^{-\Delta\E_\k/T_\k} \\
x_{k} & \textup{with probability} \quad 1 - e^{-\Delta\E_\k/T_\k} \\
\end{cases}
\end{equation*}
which depends on the temperature $T_\k$ and the gap $\Delta \E_\k := \E(y_{k+1}) - \E(x_\k)>0$. We note that different acceptance probabilities have been proposed, see \cite{romeo1991} for a discussion.
\item Before re-iterating the procedure, the system is cooled by decreasing the value of the parameter $T_\k$. A classical choice is to set  
\[T_\k = \frac C {\log(k + 2)}\]
 for some $C>0$ \cite{hajek1988}, but  stochastic or adaptive choices are also possible \cite{belisle1992}.
\end{itemize}

The intuition behind the optimization strategy is that large fluctuations given by high temperatures, allow to escape local minima. Therefore, high acceptance probabilities are set at the early stage of the computation
to favor exploration. During the computation, the acceptance rules become more strict as the system cools down. Ultimately, the algorithm will tend to accept only candidate solutions that lower the objective value. The temperature $T_\k$, therefore, is the key parameter of the dynamics, and its update rule determines the balance between exploration and exploitation in SA algorithms.

%random walk?

\begin{algorithm}
\caption{Simulated Annealing}
\label{alg:sa}
\hspace{-3mm}
\DontPrintSemicolon

$x_{0} \in \Rd, T_{0}>0, k = 0$  \tcc*[r]{\scriptsize{initialization}}
\Repeat{stop condition is satisfied}{
 $T_{k}  = 1/\log(k+2)$\tcc*[r]{\scriptsize{set decreasing temperature}}
 $y_{k+1}  = x_\k + \sqrt{2T_\k}\xi_\k$ \tcc*[r]{\scriptsize{candidate update}}
\uIf{$\E(x_{k+1}) < \E(y_\k)$}
{$x_{k+1} = y_{k+1}$ \tcc*[r]{\scriptsize{candidate solution accepted}}}
\Else{
	$\beta \sim \textup{Unif}[0,1]$\;
	\uIf{$\beta < \exp( (\E(y_{k+1}) - \E(x_\k))/T_\k) $}{$x_{k+1} = y_{k+1}$ \tcc*[r]{\scriptsize{candidate solution accepted}}}
	\Else{$x_{k+1} = x_\k$ \tcc*[r]{\scriptsize{candidate solution not accepted}}}
	}
      $k = k+1$\;
  }
  \KwOut{$x_{\k}$ \tcc*[r]{\scriptsize{return last iterate}}}
  \end{algorithm}

\subsection{A linear kinetic model}
 
In this section, we show how to formally derive a time-continuous kinetic model describing the SA dynamics. We also show how to obtain the Langevin dynamics under a suitable scaling of the parameters.

Let $f_{0}$ be the initial probability density over $\Rd$ used to sample the initial candidate solution. 
Consider the acceptance function
\begin{equation}
\beta(x,y) := \min \left \{1, e^{- (\E(y) -\E(x))/T}\right \}\,.
\end{equation}
It is immediate to verify that a single run of SA corresponds to a single realization of the Markov chain initialized with $x_{0} \sim f_{0}$ and for $k = 0, 1, \dots$
\begin{equation}
x_{k+1} = x_\k + \gamma_\k \sigma \xi_\k
\end{equation}
with 
$\gamma_\k \sim \textup{Bern}\left(\beta\left( x_\k, x_\k + \sigma \xi_k \right)\right).$
For any test function $\phi \in \mathcal{C}_b(\Rd)$, the corresponding Markov transition operator is given by
\begin{align}
Q\phi(x) &=\mathbb{E}\left[ \phi(x_{k+1})\; \middle | \;x_\k = x \right] \notag \\
& = \mathbb{E}\left[ \phi(x+ \sigma \xi_\k) \mathbb{P}(\gamma_\k = 1) + \phi(x )\mathbb{P}(\gamma_\k = 0) \right] \notag\\
& = \mathbb{E}\left[ \phi(x + \sigma \xi_\k) \beta(x, x+ \sigma \xi_\k ) + \phi(x)\left(1 -  \beta(x,x+ \sigma \xi_\k ) \right) \right] \notag \\
& =  \mathbb{E}\left[ \left( \phi(x + \sigma \Xi_\k)  - \phi(x) \right) \beta(x,x+ \sigma \Xi_\k )\right]  +\phi(x)  \notag \\
& =  \mathbb{E}\left[ \beta(x,x+\sigma \xi) \left( \phi(x+\sigma\xi) - \phi(x) \right)\right] + \phi(x) \label{eq:Qsa}
\end{align}
where the last expectation is taken with respect to $\xi\sim \mathcal{N}(0,I_d)$.

We are interested in deriving a corresponding time-continuous model $f = f(x,t)$ and parameter $\sigma = \sigma(t)$.  Assuming to expect one update per unit time, after a time step $\Delta t \in(0,1)$ the particle distribution is modified as
%\begin{equation}
%f_{t + \Delta t} = \Delta t Q^* f_t + (1 - \Delta t) f_t
%\end{equation} 
%or, equivalently,
\begin{align*}
\langle f (\cdot,t + \Delta t), \phi \rangle  &= \Delta t\langle f(\cdot,t), Q \phi \rangle + (1-\Delta t) \langle f(\cdot,t) ,  \phi \rangle \\
& = \Delta t \mathbb{E}\left [\int \beta(x,x+\sigma(t) \xi) \left( \phi(x+\sigma(t) \xi) - \phi(x) \right) f(x,t)\d x \right ] \\
& \qquad - \langle f(\cdot,t), \phi \rangle\,,
\end{align*}
with $\langle \cdot,\cdot \rangle$ being the usual $L^2(\Rd)$ product.
By formally taking $\Delta t \to 0^+$, we have
\begin{equation}
\frac{\d}{\d t} \int \phi (x) f(x,t)\d x  = \mathbb{E}\left[\int \beta(x,x') \left( \phi(x') - \phi(x) \right) f(x,t)\d x \right]
\label{eq:kinsa}
\end{equation}
where we introduced the concise notation $x' = x + \sigma(t) \xi$. By exploiting the symmetry of the random perturbation $\xi \sim \mathcal{N}(0,I_d)$, see \cite{pareschi2024sa}, \eqref{eq:kinsa} in strong form
reads
\begin{equation}
    \partial_t f (x,t) = \mathcal{L}[f](x,t)
\end{equation}
with $\mathcal{L}[f](x,t) := \mathbb{E}\left[b(x',x) f(x',t)\right]
     - \mathbb{E}\left[b(x,x') f(x,t)\right]$.

\subsection{Convergence results}
\label{sec:sa:convergence}
SA is arguably the most studied stochastic particle method in optimization, thanks to its connections with the Metropolis--Hasting sampling method. To our knowledge, the time-continuous model \eqref{eq:kinsa} was first introduced in \cite{holley1988} to study the SA for optimization over finite states $x$. In \cite{pareschi2024sa} it was extended for optimization over $\Rd$.

The analysis of the kinetic model is based on the evolution of the relative Shannon--Boltzmann entropy (also known as Kullback--Leibler divergence)
\begin{equation*}
H(g | f_\infty^T)= \int f_\infty^T(x) \log\left(\frac{g(x)}{f^T_\infty(x)} \right)\d x \,,
\end{equation*}
for a probability density $g$. In particular, given the quantity
\begin{equation}
I[g] = \frac12 \mathbb{E}\left[\iint\frac{\exp\left(-\max\{\E(x),\E(y)\}/T\right)}{Z_T} h\left(\frac{g(y)}{f_\infty^T(y)},\frac{g(x)}{f_\infty^T(x)} \right)\, \d x\, \d y \right], 
\label{eq:saI}
\end{equation}
with $h(x,y) = (x-y)(\log x - \log x) \geq 0$, for 
time-dependent temperature $T(t) = O(1/\log(t))$ and diffusion $\sigma(t) = O(\sqrt{T(t)})$ it holds:
\begin{equation}
\frac{\d}{\d t} H(f(\cdot,t) | f_\infty^T) \leq - I[f(\cdot,t)] +  \frac{C}t \| \E\|_{L^\infty}\,.
\label{eq:entropydt}
\end{equation}
Therefore, estimates of the type \cite{villani2001entropic}
\begin{equation}
I[f(\cdot,t)] \geq \lambda(t) H[f(\cdot,f)| f_\infty^T]\,,
\label{eq:salog}
\end{equation}
for some $\lambda(t) >0$, which are based on log-Sobolev inequalities, are central in the analysis of SA algorithms \cite{holley1988, pareschi2024sa}. Indeed, for a solution $f$ to \eqref{eq:kinsa}, the following convergence result holds.
\begin{theorem}[{\cite[Theorem 3.3]{pareschi2024sa}}]
If $\lambda = \lambda(t)$ satisfies \eqref{eq:salog} with $\lambda(t) = C_2 t^\alpha$, $C_2>0, \alpha \in (0,1)$, then for $T(t) = (C_1\log(t))^{-1}$, $C_1>0$, we have for $t\geq t_0$ sufficiently large
\[
H(f(\cdot,t)|f^{T(t)}_\infty) \leq \frac{2C_1}{C_2} \|\E \|_{L^\infty} t^{-\kappa} + H(f(\cdot,t_0)|f^{T(t_0)}_\infty)\exp\left( -\frac{C_2}{\kappa}(t^{\kappa} - t_0^{\kappa})\right)
\]
where $\kappa = 1-\alpha$.
\end{theorem}

\subsection{Diffusion scaling}

\label{sec:sa:limit}
Inspired by the SA mechanism, in \cite{stuart1986diffusions} the authors propose a noisy gradient decedent algorithm that follows the overdamped Langevin dynamics
\begin{equation}
\partial_t f(x,t) =  \nabla \cdot (\nabla \E(x) f(x,t)) +  \sigma(t)^2 \Delta f(x,t)\,.
\label{eq:langevin}
\end{equation}
Model \eqref{eq:langevin} describes the evolution of a particle in a viscous fluid and, as for the Metropolis--Hasting Markov chain, the Boltzmann--Gibbs density $f^T_\infty$ \eqref{eq:gibbs} corresponds to the equilibrium for $\sigma = \sqrt{2T}$ and a fixed $T>0$.

In \cite{gelfand1987analysis} it is shown that a suitable interpolation of the SA algorithm converges to \eqref{eq:langevin}. Interestingly, such a connection can also be shown at the  kinetic level. In particular, as illustrated in \cite{pareschi2024sa}, by taking the rescaling 
\[t \to t/\ve, \quad \sigma \to \sqrt{\ve}\sigma\]
and the limit $\ve \to 0^+$, the linear kinetic equation  \eqref{eq:kinsa} converges to the Fokker--Planck equation \eqref{eq:langevin}. Such scaling limit was inspired  by the so-called grazing collision asymptotics of the Boltzmann equation \cite{desvillettes1992asymptotics}.
Related diffusion limits are also known for random walk Metropolis algorithms
\cite{gelman1997weak,stuart2012diffusion}.

We will show here the heuristic behind the argument in $d = 1$ and a fixed $\sigma >0$, for illustrative purposes. Under the diffusive scaling and for $\phi\in \mathcal{C}^\infty_0(\R)$, \eqref{eq:kinsa} can be written as
\begin{equation}
\frac{\d}{\d t}\int \phi(x) f(x,t)\d x  =\frac1{\ve}
 \int \mathbb{E}\left[\phi(x+\sigma \sve\gamma \xi) - \phi(x) \right] f(x,t)\d x
\label{eq:saeps}
\end{equation}
with $\gamma \sim \textup{Bern}(x + \sve \sigma \xi)$,  $\xi \sim \mathcal{N}(0,1)$, and the expectation taken with respect to both random variables.
Consider the Taylor's expansion as $\ve \to 0^+$:
\begin{equation}
\mathbb{E} \left [\phi(x + \sve \sigma \gamma \xi) - \phi(x) \right]  = 
\sve \sigma  \mathbb{E} \left [ \gamma \xi\right]   \partial_x \phi(x) + \frac{\ve\sigma^2}2  \mathbb{E}\left[ (\gamma \xi)^2 \right] \partial_{x}^2 \phi(x) + o(\ve)\,.
\label{eq:sa1}
\end{equation}
To study the asymptotic behaviour of $\mathbb{E} \left [ \gamma \xi\right]$ and $\mathbb{E} \left [ (\gamma \xi)^2\right]$ we exploit the fact that $\min\{1,  e^z\}$ is $1$-Lipschitz. We have, therefore, 
\begin{align*}
\mathbb{E} \left [ \gamma \xi\right] &= \mathbb{E} \left [ \beta(x,x+\sve\sigma\xi) \xi\right] \\
&= \mathbb{E} \left [ \min\left\{1 , \exp(\E(x)/T - \E(x + \sve\sigma  \xi )/T\right\} \xi\right] \\
& =  \mathbb{E} \left [ \min\left\{1, \exp(- \sve\sigma \xi \partial_x \E(x)/T )\right\} \xi\right] + o(\sve)\,.
\end{align*}
From \cite[Lemma 2.4]{stuart2012diffusion}, we recall that if $X \sim \mathcal{N}(0,1)$ and $a \in \R$ it holds
\[
\mathbb{E}\left[X\min\{1,  \exp(aX)\}\right] =  a \exp(a^2/2) \Phi(-|a|)\,.
\]
with $\Phi$ being the cdf of the standard normal distribution. It follows 
\begin{align*}
\mathbb{E} \left [ \gamma \xi\right] & = -\sve \sigma \partial_x \E(x)/T \exp\left(\ve (\sigma/T)^2 (\partial_x \E(x))^2/2 \right) \Phi(- |\sve \sigma \partial_x \E(x)/T |)+ o(\sve) \\
& = -\sve \sigma \partial_x \E(x)/(2T) + o(\sve)
\end{align*}
given  $\exp(a^2/2) \Phi(- |a|) \to 1/2$ as $a \to 0$. For $\mathbb{E} \left [ (\gamma \xi)^2\right]$, since $\mathbb{E}[\gamma^2] = \mathbb{E}[\gamma]$
we have 
\begin{align*}
\mathbb{E} \left [ (\gamma \xi)^2\right] & = \mathbb{E} \left [\min\left\{1,  \exp\left(\E(x)/T - \E(x + \sve\sigma  \xi )/T\right) \right\} \xi^2\right] + o(\sve) \\
&  =  \mathbb{E} \left [\min\left\{1,  \exp\left( - \sve \sigma \xi \partial_x\E(x)/T\right) \right\} \xi^2\right]  + o(\sve)\,.
\end{align*}
We now use again that $\min\{1,  e^z\}$ is $1$-Lipschitz, to obtain:
\[
| \min\{1, \exp\left( - \sve \sigma \xi \partial_x\E(x)/T\right)\} - \min\{1, \exp(0)\} | = o(\sve)
\]
and so $ \mathbb{E} \left [ (\gamma \xi)^2\right] = \mathbb{E}[\min \left\{1, \exp(0)\right\}\xi^2] + o(\sve) = 1 + o(\sve)$.
By \eqref{eq:sa1}, 
\begin{align*}
\mathbb{E} \left [\phi(x + \sve \sigma \gamma \xi) - \phi(x) \right] 
%& =  \sve \sigma  \mathbb{E} \left [ \gamma \xi\right]   \partial_x \phi(x) + \frac{\Delta t\sigma^2}2  \mathbb{E}\left[ (\gamma \xi)^2 \right] \partial_{x}^2 \phi(x)\\
&=  \sve \sigma \left(-\sve \sigma \partial_x \E(x) /2 + o(\sve) \right)  \partial_x \phi(x) \\
& \qquad + \frac{\ve\sigma^2}2  \left(1 + o(\sve) \right) \partial_{x}^2 \phi(x)\\
& =  \frac{\ve\sigma^2}{2T} \partial_x \E(x) \partial_x \phi(x) + \frac{\ve\sigma^2}2 \partial_{x}^2 \phi(x) + o(\ve)\,.
\end{align*}

Therefore, as $\ve\to 0^+$, \eqref{eq:saeps} reads 
\begin{equation*}
    \frac{\d}{\d t}\int \phi(x) f(x,t)\d x
    \int \left(\frac{\sigma^2}{2T} \partial_x \E(x) \partial_x \phi(x) + \frac{\sigma^2}2 \partial_x^2 \phi(x)\right) f(x,t) \d x + o(1)\,.
\end{equation*}
and with choice $\sigma = \sqrt{2T}$ we obtain in the limit the overdamped Langevin dynamics \eqref{eq:langevin} in $d=1$.

\begin{remark}
Convergence to global minimizers for the annealed overdamped Langevin dynamics \eqref{eq:langevin} follows the same arguments illustrated in Section \ref{sec:sa:convergence}, but the relevant quantity 
 $I$ \eqref{eq:saI} corresponds in this case to the Fisher information \cite{monmarche2018,tang2024discrete,chizat2022meanfield}.
\end{remark}

\section{Genetic Algorithm}
\label{sec:ga}
% use word jargon
\subsection{The heuristic strategy}
\label{sec:ga:heuristics}

Genetic Algorithms (GA) \cite{holland1975} are population-based methods inspired by the mechanism of evolution and natural selection of biological organisms. In the GA terminology, every individual of the ensemble has a certain \textit{fitness} which depends on the objective value. Intuitively, a low objective value corresponds to a high fitness level. The fittest individuals are then selected for reproduction. Two \textit{parents} individuals, generate a \textit{child} individual via the biological-inspired operators \textit{crossover}, and \textit{mutation}. Crossover aims to mix the features of the two parents, while the mutation operator randomly modifies the features. We note that GA was first proposed for individuals represented by binary strings, but its computational paradigm can be extended to real-valued vectors, which is the setting we will consider here. We refer to the recent survey \cite{ga2021review} for more details on GAs.

Let us describe more in detail the GA algorithmic components. At step $k = 0,1,2 ,\dots$ the population ensemble is described by $N$ individuals $\textbf{x}_\k = (x_\k^1,\dots,x^N_\k)$, $x^i_\k\in \Rd$, $i = 1, \dots, N$, which we will call \textit{particles} from now on. The iterative evolution is then characterized by the following steps:

\begin{itemize}
\item \textbf{Selection of fittest particles.} To every particle $x^i_\k$ we assign a fitness value $\Par^i[\textbf{x}_{\k}]$ such that $\sum_{i=1}^N \Par^i[\textbf{x}_{\k}] = 1$. The choice of the fitness function determines the selection mechanism of the GA. The fitness can be attributed to the particle $x^i_\k$ proportionally to its objective value $\E(x_k^i)$, or to its rank in the ensemble. Popular choices for objective-based fitness are, for instance, the Boltzmann--Gibbs Selection and Random Wheel Selection given respectively by
\begin{equation} \Par^i[\textbf{x}_{\k}]  \propto \exp(-\alpha \E(x_\k^i))\qquad \textup{and} \quad \Par^i[\textbf{x}_{\k}]  \propto \max \E -\E(x_\k^i)\,,
\label{eq:selectionN}
\end{equation}
for some $\alpha>0$. Rank-based Selection is simply determined by $\Par^i[\bm{x}_{\k}]  \propto \rank(x_\k^i)$, where $\rank(x_\k^i) = j$ if $x_{\k}^i$ has the $j$-th highest objective value among all particles. We refer to \cite{khalid2013selection} for an extensive review of different selection mechanisms considered in the literature. 
Once the fitness values are determined, a pool of parents particles are randomly sampled with probabilities $\Par^1[\textbf{x}_{\k}],\dots, \Par^N[\textbf{x}_{\k}]$.

\item \textbf{Crossover and mutation.} Parent particles are used to generate new particles via binary interactions. The interaction is composed of the crossover operator which aim to mix the entries of the parents and a mutation operator which adds noise to the interaction.
The crossover depends on a possibly random parameter $\gamma \in [0,1]^d$, while mutation is determined by a parameter $\sigma>0$ and a random vector $\xi$ sampled from a mean-zero distribution with unit variance. Let $x$, $x_*$ be two particles selected from the pool of parents at step $k$, the newly generated particle is given by
\begin{equation}
x' = (1- \gamma) \odot x + \gamma \odot x_* + \sigma \xi\,,
\end{equation}
with $\odot$ being the component-wise (Hadamard) product.
\item \textbf{New population.} The created particles are included in the ensemble, but, differently to classical models of binary interactions, they do not substitute the parent particles. Still, to avoid the number of particles increasing at every iteration, a certain number of them is discarded: they can be chosen at random, or among the less fit particles. The latter strategy is known as the \textit{elitist} strategy as it keeps only the best particle from the previous population.

The new ensemble  $x^1_{k+1},\dots, x^N_{k+1}$ is therefore made partially of new generated particles and partially of particles from the previous population. The ratio between the two components is determined by a parameter $\nu \in [0,1]$.
\end{itemize}

At the end of the iterative computation, the particle with the lowest objective value is returned as the computed solution. Many different variants of the GA computational paradigm have been proposed in the literature, with different implementations of the above three main steps, see \cite{ga2021review}. We summarize in Algorithm \ref{alg:ga} the GA method which aim to model in the next section.

\begin{algorithm}
\caption{Genetic Algorithm}
\label{alg:ga}
\hspace{-3mm}
\DontPrintSemicolon
\begin{tabular}{ l l l }
$N \in \mathbb{N}$ & number of particles \\  
$\sigma \in (0,\infty)$ & mutation strength \\
$\nu \in (0,1)$ & fraction of new particles per iteration \\
\end{tabular}\tcc*[r]{\scriptsize{parameters}}

$x_{0}^i \in \Rd$ for $i = 1, \dots, N$, $k=0$ \tcc*[r]{\scriptsize{initialization}}

\Repeat{stop condition is satisfied}{
 
Evaluate $\E(x_\k^1),  \dots, \E(x_\k^N)$ \tcc*[r]{\scriptsize{fitness computation}}
Compute fitness weights $\Par^1[\textbf{x}_\k], \dots, \Par^N[\textbf{x}_\k]$\;
  \For{$i = 1, \dots, \lfloor \nu N\rfloor$}{
Sample $x$, $x_*$ from $(x^i_{\k})_{i=1}^N$ \tcc*[r]{\scriptsize{select parents}}
$\qquad$with probabilities $(\Par^i[\textbf{x}_\k])_{i=1}^N$ \;
$\gamma \sim \textup{Unif}([0,1]^d)$ \tcc*[r]{\scriptsize{sample crossover vector}}
$\xi \sim \mathcal{N}(0,I_d)$ \tcc*[r]{\scriptsize{sample mutation vector}}
$x^i_{k+1} = (1- \gamma) \odot x + \gamma \odot x_* + \sigma \xi$ \tcc*[r]{\scriptsize{generate new particle}}
}
Take $x_{k+1}^i$, $i = \lceil \nu N  \rceil, \dots, N$ randomly form $(x_\k^i)_{i = 1}^N$  %\tcc*[r]{\scriptsize{complete population}}

\qquad (possibly with fitness weights $(\Par^i[\textbf{x}_\k])_{i=1}^N$) \tcc*[r]{\scriptsize{elitist strategy}}
    $k = k+1$\;
  }
  \KwOut{$x^\best_{\k}  = \argmin \{\E(x)\mid x = x_\k^i, i = 1, \dots, N \}$ \tcc*[r]{\scriptsize{best particle}}}
  \end{algorithm}

%%% NOTES %%%
% it could make sense to present the algorithm in a more abstract level, so that we are more free in the modeling procedure
%population-based
%introduced in \cite{holland1975}
%The biological-inspired operators are selection, mutation, and crossover.
%GA review
% evolutionary algorithm
% (1) evaluation of individual fitness, (2) formation of a gene pool, and (3) recombination and mutation.
% originally binary string
%It is possible to have multiple copies of an individual in a population and/or gene pool.
%gaussian mutation
%uniform/arithmetical crossover, "exchange tails"
%Holland [38] originally suggested that one parent be chosen from the gene pool and that the second parent be chosen at random from the initial population.
% (a) mutations, (b) crossovers, and (c) selection pressure.
%the individuals are mostly coded as bitstrings in the GAs, whereasthey are mostly coded as a set of real numbers in the ESs

\subsection{A Boltzmann-type kinetic model}

We model the particles generated by Algorithm \ref{alg:ga} as a Markov process over $(\Rd)^N$, where $\textbf{x}_\k = (x^1_\k, \dots, x^N_\k)$ denotes the system at the $k$-th algorithmic step. Our aim is to first approximate the $N$ particle system with a time-discrete kinetic description by assuming propagation of chaos, and then derive a time-continuous kinetic model.

According the evolution described in Algorithm \ref{alg:ga}, if we assume $ \nu N \in \mathbb{N}$, for any measurable test function $\phi:\Rd \to \R$ and a randomly picked particle $i\in\{1,\dots,N\}$ it holds
\begin{multline}
\mathbb{E}\left [   \phi(x_{k+1}^i) \; \big | \; \textbf{x}_\k  = \textbf{x}    \right]   = \nu \mathbb{E}\left[ \sum_{j,\ell = 1}^N \phi(\coll(x^j, x^\ell, \gamma, \xi)) P^i[\textbf{x}] P^j[\textbf{x}] \right]  \\
+ (1-\nu) \sum_{j=1}^N \phi(x^j) P^j[\textbf{x}]\,.
\label{eq:ga:xi}
\end{multline}

We note that, above, the right-hand side only depends on $\textbf{x}$ via its correspondent weighted empirical distribution $\sum_{i=1}^N P^i[\textbf{x}] \delta_{x^i}$. To derive the GA kinetic model we generalize the parent's distribution by introducing the function $\Parp: \mathcal{P}(\Rd) \to \mathcal{P}(\Rd)$ which maps any $\mu \in \mathcal{P}(\Rd)$ to the corresponding parent's measure $\Parp[\mu]$. In the case of Boltzmann--Gibbs Selection, $\Parp[\mu]$ it is defined by 
\begin{equation}
\Parp[\mu](\d x)  =  \frac{\exp(-\alpha \E(x))}{\int\exp(-\alpha \E(y))\mu(\d y)}\mu(\d x)
\label{eq:bgselection}
\end{equation}
so that when evaluated in $f^N = (1/N)\sum_{i = 1}^N \delta_{x^i}$ we obtain exactly $\Parp[f^N] = \sum_{i=1}^N P^i[\bm{x}] \delta_{x^i}$, see \eqref{eq:selectionN}. The same can be done for Random Wheel Selection. For the Rank-based selection mechanism, we note that the rank function can be re-written as
\[\rank(i) = \#\{ x^j, i =1, \dots, N \;|\; \E( x^i) \leq \E(x^j)\} = N  \int_{\E(x^i) \leq \E(y)} f^N(\d x)\,.
\]
Therefore, the definition of parent's probability measure can be extended to any $\mathcal{\mu} \in \mathcal{P}(\Rd)$ as follows:
\begin{equation}
\Parp[\mu](\d x) = \frac{ \int_{\E(x) \leq \E(y)} \mu(\d x) }{\int_\Rd\int_{\E(z) \leq \E(y)} \mu(\d y) \mu(\d z)}  \mu(\d x)\,.
\label{eq:rankselection}
\end{equation}

We are now ready to assume propagation of chaos for the particle system $\textbf{x}_\k = (x_\k^1, \dots, x_\k^N)$. Specifically, we assume that for $N \gg 1$ if $x_0^i$ are i.i.d. sampled from a common initial distribution $f_0$, the particles becomes i.i.d. also at later time steps $k\geq 0$: $\law(x_\k) = f_\k$ for some $f_\k \in \mathcal{P}(\Rd)$. As a consequence, the joint distribution of the particles is approximated as $\law(x_\k) \approx f_\k^{\otimes N}$, and 
\begin{equation}
\langle f_\k^N, \phi\rangle = \frac 1N\sum_{i=1}^N \phi(x_\k^i) \approx \langle f_\k, \phi \rangle \,, 
\label{eq:chaosga}
\end{equation}
where for $\mu\in \mathcal{P}(\Rd)$ we use the notation $\langle \mu, \phi\rangle = \int \phi(x) \mu(\d x)$.
Note that the left-hand side in \eqref{eq:chaosga} is a random variable as $f^N_\k \in \mathcal{P}(\mathcal{P}(\Rd))$, while the right-hand side is deterministic. We will also assume that $\Parp$ is sufficiently regular such that the approximation $f^N_\k \approx f_\k$ implies $\Parp[f_\k^N] \approx \Parp[f_\k]$. We refer to the next section for a discussion on the plausibility of such assumptions.

Since the particles are assumed to be independent, we study now the evolution of the mono-particle process $\OX_\k$ with $\law(\OX_\k) = f_\k$ and initial distribution $f_{0}$. Following \eqref{eq:ga:xi} and using the approximation $\sum_{i=1}^N P^i[\textbf{x}_\k]\delta_{X_\k^i}  = \Parp[f^N_\k] \approx \Parp[f_\k]$, we obtain 
\begin{multline}
\mathbb{E}\left [ \phi(\OX_{k+1}) \right] 
  =  \nu \mathbb{E}\left[ \iint \phi(\coll(x, x_*, \gamma, \xi)) \Parp[f_\k](\d x) \Parp[f_\k](\d x_*) \right]  
\\+ (1-\nu) \int \phi(x) \Parp[f_\k](\d x) \,.
\label{eq:ga:ox}
\end{multline}
It is important to note that, differently from \eqref{eq:ga:xi}, in \eqref{eq:ga:ox} there is no need to condition the expectation with respect to the configuration of the random particle system $\textbf{x}_\k$ at time $k$, as the evolution only depends on the law $f_\k$ of the particle itself, which is deterministic. 
To write more compactly the contribution of the parent's collisional interaction, we introduce for any $g\in\mathcal{P}(\Rd)$
\[
Q_c[g] := \coll_\# (\Parp[g]\otimes\Parp[g] \otimes \law(\gamma)\otimes \law(\xi))\in\mathcal{P}(\Rd) \,.
\]
From \eqref{eq:ga:ox} we have therefore that the kinetic approximation of the particle system evolves according to 
\begin{equation}
f_{k+1} = \nu  Q_c[f_\k] + (1-\nu) \Parp[f_\k] \,.
\label{eq:ga:timediscr}
\end{equation}

We would like now to model the evolution as a time-continuous process. This can be done, for instance, by assuming that the time between updates $T_\k$ is given by an external Poisson process with rate $1$. With such a rate, it means that we expect the whole population of particles to be updated per unit time, as in Algorithm \ref{alg:ga}. The probability that the kinetic approximation $f\in \mathcal{C}([0,T],\mathcal{P}(\Rd))$ is updated between time $t$ and $t+\Delta t$ is given by $\Delta t$, for some $\Delta t \in(0,1)$, leading to 
\begin{align*}
f(t + \Delta t) &= \Delta t \big( \nu  Q_c[f(t)] + (1-\nu) \Parp[f(t)] \big) + (1- \Delta t) f(t) \\
& = \Delta t \big( \nu Q_c[f(t)] +  (1- \nu) \Parp[f(t)] -  f(t) \big) + f(t)\,.
\end{align*}
By formally taking the limit $\Delta t \to 0^+$, we obtain the time-continuous kinetic model for GA algorithm (Algorithm \ref{alg:ga})
\begin{equation}
\partial_t f = \nu Q_c[f] + (1-\nu) \Parp[f] - f \,,
\label{eq:kinga}
\end{equation}
%which in weak form reads for any measurable test function $\phi$
%\begin{multline}
%\label{eq:kinga}
%\frac{\d}{\d t} \langle \phi , f(t) \rangle = \nu \mathbb{E} \left [ 
%\int_{\Rd} \int_{\Rd} \phi \left( \coll(x,x_*, \Gamma, \Xi)\right) \Parp[f(t)](\d x) \Parp[f(t)](\d x_*) \right] \\
%+ (1 - \nu) \langle \Parp[f(t)], \phi \rangle - \langle f(t) ,  \phi\rangle .
%\end{multline}

To better draw a connection with respect to the classical Boltzmann equation, let consider (with a slight abuse of notation) the kinetic distribution to be a probability density $f = f(x,t)$ and let $\mathcal{B}[g]$ be the density of $\Parp[g]^{\otimes 2}$ with respect to $g^{\otimes 2}$, that is, for instance, 
\[
\mathcal{B}[g](x,x^*)= \frac{\exp\left( - \alpha (\E(x) + \E(x^*)) \right)}{\left(\int \exp(\alpha \E(x)) g(x)\d x\right)^2}
\]
in the case of Boltzmann--Gibbs selection. Then, \eqref{eq:kinga} reads in weak form for $f = f(x,t)$
\begin{multline}
    \frac{\d}{\d t}\int \phi(x) f(x,t) \d x = \\
    \int B[f(\cdot,t)](x,x^*) \big(\nu \phi(x') + (1-\nu) \phi(x) \big ) f(x,t) f(x_*,t)\d x \d x_* \\
    - \int \phi(x) f(x,t)\d x\,.
    \label{eq:kingaweak}
\end{multline}

The above  modelling procedure was proposed in \cite{borghi2024ga}, and here we have also included the elitist strategy in the derivation. The model \eqref{eq:kingaweak}, therefore, includes in its evolution all the most important genetic mechanisms outlined in Section \ref{sec:ga:heuristics}: selection of the fittest, crossover, mutation, and elitist strategy.

\begin{remark}
\label{rmk:antielite}
In \eqref{eq:kingaweak} there is a difference between the colliding particles and those that are discarded. In classical rarefied dynamics, the newly generated post-collisional particles substitute the pre-collisional ones. Here, instead, not only the colliding parents are not discarded from the ensemble, but they are also more likely to be passed into the new generation when the elitist strategy is used.

The symmetric case would correspond to the algorithm where children replace parents after the interaction. Using the same arguments as above, it is possible to see that this different choice leads to a Boltzmann-like model 
\[
 \frac{\d}{\d t}\int \phi(x) f(x,t) \d x = 
    \nu\int B[f(\cdot,t)](x,x^*) \left( \phi(x') - \phi(x) \right ) f(x,t) f(x_*,t)\d x \d x_* \,.
\]
\end{remark}

%\begin{remark}
%\label{rmk:elite}
%When the elitist strategy is not used, the population at step $k+1$ is composed of a fraction $\nu\in(0,1]$ of newly generated particles (as in the elitist case) and a fraction $1-\nu$ of particles uniformly sampled from the population at step $k$. Following the same procedure as above, we obtain the time-discrete approximation 
%\begin{equation*}
%f_{k+1} = \nu Q_c[f_\k] + (1-\nu)f_{k}
%\end{equation*}
%which leads to the time-continuous kinetic model
%\begin{equation}
%\partial_t f_t = \nu \left( Q_c[f_t]  - f_t \right)\,.
%\end{equation}
%Therefore, in the non-elitistic case, $\nu>0$ plays the role of a time reparemetrization coefficient.
%\end{remark}

\begin{remark} 
\label{rmk:albi}
A different strategy to model the pool of parents has been proposed in \cite{albi2023kinetic}. The parent's particles, instead of being sampled from $\Parp[f(\cdot,t)]$, are modelled as a separate species, the "leaders", while the remaining particles are the "followers". Let $f^L = f^L(x,t)$, $f^F = f^F(x,t)$ denote, respectively, the leaders' and followers' distributions. Their evolution is composed by a transition operator $\mathcal{T}$ describing how particles are included or excluded from the pool of parents, and by a binary collisional operator $Q$ describing the generation of children. The proposed followers-leaders dynamics therefore reads
\begin{equation}
\begin{cases}
\partial_t f^L(x,t) - \mathcal{T}[f^L](x,t) = Q(f^L, f^F)(x,t) + Q(f^L, f^L)(x,t)  \\
\partial_t f^F(x,t) - \mathcal{T}[f^F](x,t) = Q(f^F, f^L)(x,t) + Q(f_t^F, f^F)(x,t)  
\end{cases}
\label{eq:albi}
\end{equation}
see  \cite{albi2023kinetic} for more details.
\end{remark}

\subsection{Discussion on propagation of chaos assumption}

Numerical evidence presented in \cite{borghi2024ga} suggests that the GA particle dynamics converge to the kinetic model \eqref{eq:kinga} as $N \to \infty$. Yet, whether a propagation of chaos property of type \eqref{eq:chaosga} holds remains an open question.

A main challenge in the analysis lays in the nonlinearity of the collisional kernel. The collisional frequency of a particle $x$ depends, indeed, nonlinearly on $f(t)$ via the normalization constant used in the definition of $\Parp[f(t)]$ (see \eqref{eq:bgselection} and \eqref{eq:rankselection}). Kinetic equations with such a feature were considered in \cite{degond2013kinetic} to study biological swarm models as $N \to \infty$ through the hierarchy of marginal equations. 
For such models, the authors were able to show the propagation of chaos
 for pair selection probabilities satisfying
\[
p^{i,j} \leq \frac{C}{N(N-1)} \qquad \textup{for some} \quad C>0\,,
\]
see Remark 3.2 in \cite{degond2013kinetic}. In the GA dynamics considered, we have $p^{i,j} = \Par^i[\textbf{x}]\Par^j[\textbf{x}]$. We note that, for bounded objective function $\E$ and Boltzmann-Gibbs selection, it holds for all $i = 1, \dots,N$
\[\Par^i[\textbf{x}] \leq \frac{\exp(\alpha(\max \E- \min E))}N
\]
and so the marginals-based technique proposed in \cite{degond2013kinetic} seems to be promising to prove the propagation of chaos for the GA algorithm, at least under such settings. 

We conjecture that another viable approach consists of  coupling the GA particle system $(x_\k^1,\dots,x_\k^N)$ with $N$ independent copies of the mono-particle process \eqref{eq:ga:ox} $(\OX_\k^1, \dots, \OX_\k^N)$, as outlined, for instance, in \cite{fournier2016nanbu}. We refer to the reviews \cite{diez2022review1,diez2022review2} for more details on the rigorous derivation of chaos properties for classical and non-classical particle systems. 

We note that for modelling GA algorithms we have considered the large particle limit $N \to \infty$, and, only afterward, taken the time-continuous model by considering $\Delta t \to 0^+$. The illustrated approach of linking the algorithm and kinetic model is related to Monte Carlo approximation method, and in particular the Nanbu approach \cite{nanbu1980}, where the particle method aims to approximate the dynamics via a probabilistic interpretation of the time-discretized Boltzmann equation \cite{wagner2005, pareschi13}. As already mentioned, an advantage of working with the kinetic model is that we can easily derive novel algorithms by using a different numerical approximation of the model. In the case of GA, for instance, a novel algorithm would consist of deriving first a (time-continuous) particle approximation to \eqref{eq:kinga} which is then discretized in time. 
In this case, the stochastic particle system is described by a $(\Rd)^N$-valued Markov jump process $\bm{X}_t= (X_t^1, \dots, X_t^N)$  which is numerically approximated by sampling from Poisson random measures, see \cite[Chapter 12]{del2017stochastic} and next remark for more details

\begin{remark}
\label{rmk:sdega}
A time-continuous particle system that approximates \eqref{eq:kinga} can be defined as follows. Consider infinitesimal generator $\mathcal{L}_N$ defined, for $\phi:(\Rd)^N \to \R$ sufficiently regular, as
\textup{
\begin{multline}
\mathcal{L}_N \phi(\textbf{x}) 
 = \nu   \mathbb{E}\Bigg[ \sum_{i,j,\ell=1}^N \left( \phi (\textbf{x} + (\coll(x^j,x^\ell, \gamma, \xi) - x^i)\textbf{e}_i) - \phi(\textbf{x}) \right)
\Par^{j}[\textbf{x}]\Par^{\ell}[\textbf{x}] 
\Bigg] \\
\qquad +  (1-\nu)\sum_{i,j=1}^N \left( \phi (\textbf{x} + (x^j - x^i)\textbf{e}_i) - \phi(\textbf{x}) \right)
\end{multline}}
where for $h\in \Rd$ we set \textup{$h\textbf{e}_i = (0, \dots,0, h, 0, \dots, 0) \in (\Rd)^N$} the vector with $h$ in the $i$-th place. The law $F(t) = \law(\bm{X}_t) \in \mathcal{P}((\Rd)^N)$ of the associated jump Markov process satisfies
\begin{equation*}
\frac{\d}{\d t}\langle F(t), \phi \rangle = \langle F(t), \mathcal{L}_N \phi  \rangle  \,.
\label{eq:sdega}
\end{equation*}
By taking as a test function \textup{$\phi(\textbf{x}) = \phi(x^1)$} and assuming propagation of chaos in the sense of \eqref{eq:chaosga} ($F(t) \approx f(t)^{\otimes N}$ and $\Parp[f^N(t)] \approx \Parp[f(t)]$ for $N \gg 1$), it is possible to verify that \eqref{eq:sdega} reduces exactly to the kinetic model \eqref{eq:kinga}.
\end{remark}

\subsection{Convergence result}

Convergence properties to global minimizers for the GA method described in Algorithm \ref{alg:ga} can be inferred by studying the long--time behaviour of the corresponding kinetic approximation.

In \cite{borghi2024ga}, the authors consider the time-discrete kinetic model \eqref{eq:ga:timediscr} in the case with Boltzmann--Gibbs selection of the parents and without elitist strategy. This corresponds to 
\begin{equation}
f_{k+1} = \nu Q_c[f_\k] + (1-\nu)f_{k}
\label{eq:ga:conv}
\end{equation}
where we note that, here, $\nu\in(0,1]$ can also be interpreted as a time-discretization.
The assumptions on the objective function are similar to the ones used for the analysis of Consensus-Based Optimization methods  in \cite{fornasier2024consensusbased}, namely:
\begin{assumption} \label{asm:E}
The objective function $\E: \Rd \to \R$ is continuous and satisfies: 
\begin{enumerate}[i)]
\item  \label{asm:unique} (solution uniqueness) there exists a unique global minimizer $x^\star$;
\item  \label{asm:growth}
(growth conditions) there exists $L_\E$, $c_u, c_l, R_l>0$ such that 
\begin{equation*}
\begin{cases}
|\E(x) - \E(y) | \leq L_\E (1+ |x| + |y|)|x - y| & \quad \forall\; x,y \in\Rd \\
\E(x) -  \E (x^\star)  \leq c_u (1 + |x|^2) & \quad \forall\; x \in\Rd \\
\E(x) - \E (x^\star) \geq c_l |x|^2 & \quad \forall\; x \,:\, |x|>R_l \,;
\end{cases}
\end{equation*}
\item  \label{asm:inverse}
(inverse continuity) there exist $c_p , p>0, R_p>0$  and lower bound $\E_\infty>0$ such that
\begin{equation*}
\begin{cases}
c_p|x - x^\star|^{p} \leq \E(x) - \E(x^\star) & \quad \forall x \,:\, |x-x^\star| \leq R_p \\
\E_\infty < \E(x) - \E(x^\star) & \quad \forall x \,:\, |x-x^\star| > R_p\,.
\end{cases}
\end{equation*}
\end{enumerate}
\end{assumption}
Convergence towards the minimizer is established for the mean 
\[\m[f_k] = \int x f_\k(\d x)\]
of the distribution:
\begin{theorem}[{\cite[Theorem 4.1]{borghi2024ga}}] Let $\E$ satisfy Assumption \ref{asm:E}. Let $f_{0} \in \mathcal{P}_2(\Rd)$ be such that $x^\star \in \textup{supp}(f_{0})$, and $f_{k}$ be updated according to \eqref{eq:ga:conv}. 
Fix an arbitrary accuracy $\eps>0$ and let $T^\star$ be the time horizon given by
\begin{equation*}
T^\star := \log\left(\frac{2|\m[f_\k] - x^\star|}{\eps}\right)\,.
\end{equation*}
Then, there exists $\alpha>0$ sufficiently large (the parameter in Boltzmann--Gibbs selection \eqref{eq:bgselection}) such that
\begin{equation*}
\min_{k: k\nu \leq T^\star}\left | \m[f_\k]  - x^\star \right | \leq \eps\,.
\end{equation*}
Furthermore, until the accuracy is reached, it holds
$|\m[f_\k]  - x^\star | \leq e^{-k\nu} | \m[f_0]  - x^\star|$.
\end{theorem}

The key tool for the convergence analysis is a quantitative version of the Laplace principle derived in \cite{fornasier2024consensusbased} for the study of consensus-based optimization algorithms. We conjecture that the proof strategy outline in \cite{borghi2024ga} can be extended to the study of the time-continuous kinetic equation 
\eqref{eq:kinga} with also elitist strategy.

\subsection{Quasi-invariant scaling}

Following the computations in \cite{borghi2024ga}, we show how it is possible to derive CBO-type dynamics via the so called quasi-invariant scaling \cite{pareschi13}
\[
t \to t/\varepsilon,\quad \gamma \to \varepsilon \gamma,\quad \sigma^2 \to \varepsilon \sigma^2\,,
\]
with parameter $\ve\ll1$. The scaling introduced above corresponds to a dynamic in which crossovers and mutations decrease in intensity but increase in number. 

Similarly to Section \ref{sec:sa:limit}, we consider $\phi \in C_{0}^\infty(\Rd)$ and the expansion 
\[
\phi(x')=\phi(x) + (x'-x) \cdot \nabla_x \phi(x) + \frac12(x'-x)\cdot\nabla^2\phi (x)(x'-x)+O(|x'-x|^3)\,.
\]
After the scaling, it holds
\begin{equation}
x' - x = \ve \gamma \odot (x_* - x) + \sqrt{\ve}\sigma \xi\,.
\end{equation}
By considering the case without elitist strategy ($\nu = 0$), from \eqref{eq:kingaweak} we obtain
\begin{multline*}
\frac{\d}{\d t} \int \phi (x)\,f(x,t)\d x  = \frac1{\ve} \iint \left( \mathcal{B}[f(\cdot,t)](x,x_*) - 1 \right) \phi(x)\,f(x,t)f(x_*,t)\d x \d x_* \\
+   \iint  \mathcal{B}[f(\cdot,t)](x,x_*) \gamma \odot (x_* - x)\cdot \nabla\phi(x)  \,f(x,t)f(x_*,t)\d x \d x_*\\
+  \frac12 \mathbb{E}\Big[\iiint  \mathcal{B}[f(\cdot,t)](x,x_*)  \sigma^2 \xi^\top \nabla^2\phi(x) \xi  \,f(x,t)f(x_*,t)\d x \d x_*\, \Big ] + O(\ve)\,.
\end{multline*}
To be able to formally pass to the limit $\ve \to 0^+$, we shall assume that it holds
\begin{equation*}
\int \left(\mathcal{B}[f(\cdot,t)](x,x_*) - 1 \right) \phi(x) \,f(x,t) f(x_*,t)\d x\d x_* =1\,,
\end{equation*}
which can be achieved only by considering a reproduction kernel independent on $x$,  
\[\mathcal{B}[f(\cdot,t)](x,x_*) = \mathcal{B}[f(\cdot,t)](x_*) = 
\frac{\exp(-\alpha \E(x_*))}{\int \exp(-\alpha \E(x))f(x_*,t)\d x_*}
\,.\]
This is equivalent to the situation in which only one of the two parents $(x,x_*)$ is chosen by the Boltzmann--Gibbs selection mechanism, while the other $x$ is chosen at random.

By passing to the limit $\ve \to 0^+$ and reverting to the strong form we get the Fokker--Planck equation 
\begin{equation*}
\partial_t f(x,t) 
+ \nabla_x \cdot \left(\mathcal{R}[f(\cdot,t)](x) f(x,t)\right) = \frac{\sigma^2}{2}\Delta f(x,t)
\label{eq:mf1}
\end{equation*}
with
\begin{equation*}
\mathcal{R}[f(\cdot,t)](x) = \gamma \odot \left(\frac{\int  x_* \exp(-\alpha \E(x_*))f(x_*,t)\d x_*}{\int \exp(-\alpha \E(x))f(x_*,t)\d x_*} - x \right) 
   \,.
\label{eq:mf2}
\end{equation*}
Under the choice $\gamma = (\lambda, \dots, \lambda)^\top$ for some $\lambda>0$, the resulting mean-field equation \eqref{eq:mf1} models a CBO-type evolution with non-degenerate diffusion.

%%% PSO
\section{Particle Swarm Optimization}
\label{sec:pso}
\subsection{The heuristic strategy}

The Particle Swarm Optimization (PSO) algorithm has been introduced in  \cite{kennedy1995,kennedy1997} with the objective of designing a simple computational paradigm inspired by social models. The key idea is to consider the optimization problem as a collective task between individuals (the particles) that influence each other to reach a common goal. The toy model proposed in \cite{kennedy1995}, in particular, draws inspiration from a flock in which the birds align their velocity towards food sources. 

The method considers $N$ particles which update their states at every iteration $k = 1, 2, \dots $ of the algorithm. With $x_\k^i \in \Rd$ with denote the position of the $i$-th particle, and with $v^i_\k \in \Rd$ its velocity. 
The update of the velocity $v^i_\k $is characterized by the following elements:
\begin{itemize}
\item Alignment towards the best position (with respect to the objective $\E$) found by the particle $i$ until the present step $k$. In the PSO nomenclature, this is called the \textit{personal best} and we indicate it with $y^i_{\k} \in \Rd$.
\item Alignment towards the \textit{global best} position found by the whole ensemble until step $k$, which we indicate with $y^\best_\k \in \Rd$.
\item Two random weights $r^{i,1}_\k, r^{i,2}_\k \sim \textup{Unif}([0,1]^d)$ and two parameters $c_1, c_2 >0$ called \textit{acceleration coefficients}.
\end{itemize}
The update of the particle position is characterized by a simple free transport evolution. Altogether, the PSO update rule reads
\begin{equation}
\begin{cases}
x_{k+1}^i &= x_{\k}^i  + v^i_{k+1} \\
v_{k+1}^i & = v_\k^i + c_1 r^{i,1}_\k \odot (y^i_\k - x^i_\k) + c_2 r^{i,2}_\k\odot (y^\best_\k - x_k^i)
\end{cases}
\qquad i = 1, \dots,N\,.
\label{eq:pso}
\end{equation}
We notice that the new velocity $v_{k+1}^i$ is given by a combination between the previous velocity $v_\k^i$, and the two differences $(y^i_\k - x^i_\k)$, $(y^\best_\k - x_k^i)$, which are also multiplied component-wise with the random vectors $r^{i,1}_\k, r^{i,2}_\k$.  The suggested choice for  the acceleration coefficients is $c_1= c_2 = 2$, so that the expected value of $c_1r^{i,1}_\k$, $c_2r^{i,2}_\k$ is the unit vector
\begin{equation}
\mathbb{E}[c_1 r_{k}^{i,1}] = \mathbb{E}[c_2 r_{k}^{i,2}] = (1, \dots, 1)^\top\,.
\label{eq:pso_c1r1}
\end{equation}
Let $H[\cdot]$ denote the Heaviside step function ($H[z] = \bm{1}_{[0,\infty]}(z)$), and the update of the personal bests can be implemented as follows:
\begin{equation}
y^i_{k+1} = y_\k^i + H\left[ \E(y_\k^i) - \E(x_{k+1}^i)\right] (x_\k^i - y_\k^i)\,.
\label{eq:pso_pbest}
\end{equation}
Since its introduction, several different variants and improvements of PSO have been designed: we refer to the recent review \cite{pso2022survey} for a detailed overview of the topic. We summarize the fundamental PSO strategy in Algorithm \ref{alg:pso}.

\begin{algorithm}
\caption{Particle Swarm Optimization}
\label{alg:pso}
\hspace{-3mm}
\DontPrintSemicolon
\begin{tabular}{ l l l }
$N \in \mathbb{N}$ & number of particles \\  
$c_1 \in (0,\infty)$ & personal best acceleration coefficient \\
$c_2 \in (0,\infty)$ & global best acceleration coefficient \\
\end{tabular}\tcc*[r]{\scriptsize{parameters}}

$x_{0}^i, v_{0}^i, y_{0}^i \in \Rd$ for $i = 1, \dots, N$ \tcc*[r]{\scriptsize{initialization}}

$y^\best_{0} =  \argmin \{ \E(y) \mid y = y^i_{0}, i =1, \dots, N \}$ 

$k = 0$\;

\Repeat{stop condition is satisfied}{
 
  \For{$i = 1, \dots,N$}{
$r^{i,1}_\k, r^{i,2}_\k \sim \textup{Unif}([0,1]^d)$  \tcc*[r]{\scriptsize{sample random vectors}}
  $v_{k+1}^i  =  v_\k^i + c_1 r^{i,1}_\k \odot (y^i_\k - x^i_\k) + c_2 r^{i,2}_\k\odot (y^\best_\k - x_k^i)$\;  \tcc*[r]{\scriptsize{update velocity}} 
$x_{k+1}^i = x_{\k}^i  + v^i_{k+1}$ \tcc*[r]{\scriptsize{update position}} 
\lIf*{$\E(x_{k+1}^i) < \E(y^i_\k)$}{$y^i_{k+1} = x_{k+1}^i$} \tcc*[r]{\scriptsize{update personal best}}
\lElse{$y^i_{k+1} = y_{\k}^i$} }
  $y^\best_{k+1} = \argmin \{\E(y)\mid y = y_{k+1}^i, i = 1, \dots, N \}$  \tcc*[r]{\scriptsize{update global best}}
    $k = k+1$\;
  }
  \KwOut{$y^\best_{\k}$ \tcc*[r]{\scriptsize{return global best}}}
  \end{algorithm}

\subsection{Derivation of the kinetic model}
\label{sec:pso:derivation}

The objective of this section is to derive a kinetic model that approximates the particle dynamics described by \eqref{eq:pso}. We obtain this by starting from a continuous-in-time stochastic particle system and then obtain a mono-particle process under the propagation of chaos assumption.

We model the time-discrete particle system as a Markov process where each particle is described by the tuple $(X_\k^i, V_\k^i, Y_\k^i)$ taking values in $(\Rd)^3$. The update rules are modified by taking into  account the following approximations:
\begin{itemize}
\item We approximate the random variables $r_\k^{i,1}$, $r_\k^{i,2}$ as
\begin{equation}
 r_\k^{i,1} \approx  1 + \frac{1}{\sqrt{3}}\xi_\k^{i, 1} \quad \textup{and} \quad r_\k^{i,2} \approx 1 + \frac{1}{\sqrt{3}}\xi_\k^{i, 2}
\end{equation}
with  $\xi_\k^{i,1}, \xi_\k^{i,1} \sim \mathcal{N}(0, I_d)$. Such approximation preserves mean and variance of the random variable and allows to split velocity displacement in deterministic and zero-mean components.
\item We introduce a time discretization parameter $\Delta t\in (0,1]$. We rescale the deterministic displacements by $\Delta t$, and the stochastic ones by $\sqrt{\Delta t}$. As we will see, this will allow us to preserve both components in the limit $\Delta t \to 0$.
\item 
We approximate the Heaviside function used in the update of the personal best, see \eqref{eq:pso_pbest}, as
\begin{equation}
H(z) \approx  H^\beta(z) := 1 +\frac12\tanh(\beta z)\,,
\end{equation}
in order to have $H^\beta \to H$ pointwise as $\beta \to \infty$, but $H^\beta \in C^\infty(\R)$. We also add an additional parameter $\nu>0$ to control the speed of the (approximated) personal bests.
\item Finally, the global best is regularized by considering a weighted average of the personal bests. For any $\rho \in \mathcal{P}(\Rd)$ and $\alpha>0$, we introduce the average with Boltzmann--Gibbs weights
\begin{equation}
\m^\alpha_\E[\rho] :=  \frac{\int x \exp(-\alpha \E(x)) \rho(\d x)}
{\int \exp(-\alpha \E(x)) \rho(\d x)}\,.
\end{equation}
Let $\rho^{N,y}_{\k} = (1/N) \sum \delta_{Y^i_{\k}}$ be the empirical measure associated with the personal bests, we consider
\begin{equation}
y^\best_\k \approx \m^\alpha_\E[\rho^{N,y}_\k]\,.
\end{equation}
As for the personal bests, we note the approximation is exact in the limit $\alpha \to \infty$. Indeed, provided there is a unique minimum among $Y_k^1, \dots, Y_\k^N$, direct computations lead to
\begin{equation}
\m^\alpha_\E [\rho^{N,y}_{\k}] \longrightarrow \argmin \{\E(y)\mid y = Y_\k^i, i = 1, \dots, N \} \qquad \textup{as} \;\; \alpha \to \infty\,,
\end{equation}
where the above limit is to be intended event-wise, as the objects involved are random variables.
\end{itemize}

After taking into account the above approximations, the PSO update rule \eqref{eq:pso} results in 
\begin{equation}
\begin{cases}
X^i_{k+1} &= X_{k}^i + \Delta t V^i_{k+1} \\
V^i_{k+1} &  = V^i_\k  +  c_1 (Y^i_\k - X_\k^i) \odot  \left(\Delta t + \frac{\sqrt{\Delta t}}{\sqrt{3}} \xi_\k^{i,1}\right)\\
& \qquad + c_2  (\m_\E^\alpha[\rho^{N,y}_\k] - X_\k^i  )\odot \left(\Delta t + \frac{\sqrt{\Delta t}}{\sqrt{3}} \xi_\k^{i,2}\right) \\ 
Y^i_{k+1} &= Y_\k + \nu \Delta t H^\beta\left[\E(Y_\k^i) - \E(X_{k+1}^i) \right](X_{k+1}^i - Y_\k^i)
\end{cases} \qquad i = 1, \dots,N\,.
\label{eq:pso_em}
\end{equation}
Thanks to the introduction of the time step $\Delta t$, the above system can be interpreted as a semi-implicit Euler--Maruyama discretization \cite{desmond2001} of an underlying time-continuous particle system. The deterministic updates of $X^i_\k$ and $Y_\k^i$ are implicit in $V_\k^i$ and $X_\k^i$ respectively, while the update for $V_\k^i$ is fully explicit.

Let $(X_t^i, V_t^i, Y_t^i)$ for $i =1, \dots,N$ be the particle system for $t\geq 0$ and $(B_t^{i,1})_{t\geq 0}$, $(B_t^{i,2})_{t\geq 0},$ $i= 1, \dots,N$ be $2N$ independent Brownian processes. The scheme \eqref{eq:pso_em} discretizes the system of SDEs 
\begin{equation}
\begin{cases}
\d X_t^i &= V_t^i \d t \\
 \d V_t^i & = c_1 (Y_t^i - X_t^i) \d t + \frac{c_1}{\sqrt{3}} (Y_t^i - X_t^i) \odot \d B_t^{i,1}  \\
& \quad  + c_2 (\m_\E^\alpha[\rho^{N,y}(t)] - X_t^i )\d t+  \frac{c_2}{\sqrt{3}} (\m_\E^\alpha[\rho^{N,y}(t)] - X_t^i ) \odot \d B_t^{i,2} \\
\d Y^i_t &=\nu H^\beta\left[\E(Y_t^i) - \E(X_t^i) \right](X_t^i - Y_t^i)\d t
\end{cases} \quad i = 1, \dots,N\,.
\label{eq:pso_sde}
\end{equation}

We note that the communication between the particles happens uniquely in the computation of the regularized global best position $\m_\E^\alpha[\rho^{N,y}(t)]$. To study the dynamics of a single, typical, particle of the system we average the influence of the system by taking an approximation of mean-field type. 
We assume, therefore, that for large population sizes $N \gg1$ the chaos propagates from the initial data to  $t>0$: if $\law(X_0^i, V_0^i, Y_0^i) = f_0 \in \mathcal{P}((\Rd)^3)$ for $i = 1, \dots,N$  then
\begin{equation}
\law\left( (X_t^1, V_t^1, Y_t^1), \dots, (X_t^N, V_t^N, Y_t^N)\right) \approx f(t)^{\otimes N} 
\qquad \textup{for some} \quad f(t) \in \mathcal{P}((\Rd)^3)\,.
\label{eq:psochaos}
\end{equation}
Also, we obtain $\m_\E^\alpha[\rho^{N,y}(t)] \approx \m_\E^\alpha[\rho^{y}(t)]$ where $\rho^y(t)$ is the marginal of $f(t)$ over the personal best space.  Under such mean-field approximation, the dynamic of each particle depends solely on its state and on the mean-field law $f(t)$. 

System \eqref{eq:pso_sde} therefore reduces to $N$ copies of the mono-particle process of McKean type
\begin{equation}
\begin{cases}
\d\overline{X}_t &= \OV_t \d t \\
\d \OV_t^i & = c_1 (\OY_t -\overline{X}_t) \d t + \frac{c_1}{\sqrt{3}} (\OY_t - \overline{X}_t) \odot \d \OB_t^{1}  \\
& \qquad  + c_2 (\m_\E^\alpha[\rho^{y}(t)] -\overline{X}_t )\d t+  \frac{c_2}{\sqrt{3}} (\m_\E^\alpha[\rho^{y}(t)] - \overline{X}_t ) \odot \d \OB_t^{2} \\
\d\OY_t &= \nu H^\beta\left[\E(\OY_t) - \E(\overline{X}_t) \right](\overline{X}_t - \OY_t)\d t \\
 \rho^y(t) & = \law(\OY_t)
\end{cases}
\label{eq:pso_mono}
\end{equation}
where $(\OB_t^1)_{t\geq 0}, (\OB_t^2)_{t\geq 0}$ are two independent Brownian processes.

The corresponding PDE formulation is obtained by applying It{\^o}--Doeblin formula and by taking zero-value of the stochastic integrals. In strong form, the resulting kinetic model $f \in \mathcal{C}([0,T],\mathcal{P}((\Rd)^3))$ for the PSO algorithm (Algorithm \ref{eq:pso}) is given  by the non-linear Vlasov--Fokker--Planck equation
\begin{multline}
\partial_t f + v\cdot \nabla_x f   + \nabla_y \cdot \left (\nu(x - y)H^\beta(\E(y) - \E(x)) f\right )= \\ 
\nabla_x \cdot \left( c_1 (y - x) f + c_2 (\m^\alpha_\E[\rho_t^y] - x) f \right) \\
+ \frac16 \sum_{\ell=1}^d   \partial_{v_\ell}^2 \left( c_1^2(y - x)_\ell^2  f + c_2^2(\m^\alpha_\E[\rho_t^y] - x)_\ell^2  f \right)\,.
\label{eq:pso:strong}
\end{multline}

\subsection{Mean-field limit results}

 A continuous description of PSO dynamics was first derived in \cite{emara2004cpso}, while the mean-field approximation procedure just described was proposed in \cite{grassi2021pso} (with the additional inertia mechanism, see Remark \ref{rmk:inertia}). Well-posedness of the SDE particle system \eqref{eq:pso_sde} 
and of the mean-field PSO \eqref{eq:pso:strong} without personal bests have been established in \cite{huang2021psomf} assuming Lipschitz continuity and boundedness of the objective function $\E$. The existence of strong solutions for the particle systems follows from standard results in SDE theory \cite{durrett2018}, while the existence of solutions for the McKean mono-particle process is proved by applying Leray--Schauder fixed point theorem. The well-posedness results have been extended in \cite{grassi2021psosurvey} to globally Lipschitz objectives $\E$ which have a quadratic lower bound, and that grows at most quadratically at infinity. The case with also local bests have been considered in \cite{huang2023pso}.
 
In \cite{huang2021psomf}, the author proves propagation of chaos for the PSO particle system without local best update, but the proof can be directly extended to also include the memory mechanism. This is done by coupling a particle system $(X^i_t, V_t^i, Y_t^i)_{i=1}^N$ \eqref{eq:pso_sde} with $N$ copies $(\overline{X}^i_t,\OV_t^i, \OY_t^i)_{i=1}^N$ of the corresponding mean-field mono-particle process \eqref{eq:pso_mono}, and by showing that
\begin{equation}
\max_{i=1, \dots, N} \mathbb{E}\left[ |X_t^i - \overline{X}_t^i |^2 + |V_t^i - \OV_t^i|^2 +  |Y_t^i - \OY_t^i |^2 \right] \leq C N^{-1}
\label{eq:pso:coupling}
\end{equation}
within a given time window $t \in [0,T]$. The connection with the propagation of chaos assumption \eqref{eq:psochaos} lies on the fact that from \eqref{eq:pso:coupling} we obtain convergence in 2-Wasserstein metric
of the empirical measure $f^N$ of the particle system \eqref{eq:pso_sde} towards the solution $f$ to \eqref{eq:pso:strong}: $W_2(f^N(t), f(t)) \to 0$ (provided $f(t) \in \mathcal{P}_2((\Rd)^3)$). In turn, this implies  convergence in the weak sense
\begin{equation}\label{eq:convergenceW}
f^N(t) \rightharpoonup f(t) \qquad \textup{as} \quad N \to \infty\,,
\end{equation}
for $t\in [0,T]$. We refer to \cite{santa2015} for more details on Wasserstein metrics. Quantitative mean-field approximation errors are also provided in \cite{huang2023pso} for a larger class of objective functions  $\E$.

\subsection{Convergence result}

In \cite{grassi2021pso} the authors derive a more general version of the PSO dynamics where the deterministic and the stochastic components of the dynamics \eqref{eq:pso_sde} are tuned with different parameters $\lambda_1, \lambda_2, \sigma_1, \sigma_2>0$. This is conjectured to allow for a better control of the balance between exploration and exploitation with respect to the the inertia tuning mechanism. 

In \cite{huang2023pso}, the authors study the long--time behaviour of such mean-field generalized model PSO model with inertia $m>0$ and friction coefficient $\gamma = 1-m\geq0$, given by
\begin{equation}
\begin{cases}
\d\overline{X}_t &= \OV_t \d t \\
m\d \OV_t^i & = - \gamma \OV_t + \lambda_1 (\OY_t -\overline{X}_t) \d t + \sigma_1 (\OY_t - \overline{X}_t) \odot \d \OB_t^{1}  \\
& \qquad  + \lambda_2 (\m_\E^\alpha[\rho^{y}(t)] -\overline{X}_t )\d t+  \sigma_2 (\m_\E^\alpha[\rho^{y}(t)] - \overline{X}_t ) \odot \d \OB_t^{2} \\
\d\OY_t &= \nu H^\beta\left[\E(\OY_t) - \E(\overline{X}_t) \right](\overline{X}_t - \OY_t)\d t \\
 \rho^y(t) & = \law(\OY_t)
\end{cases}
\label{eq:pso_mono_gen}
\end{equation}

The emergence of consensus is done by studying the evolution of 
%\small{
\begin{multline*}
\mathcal{H}(t):=\left(\frac\gamma{2m} \right)^2|\overline{X}_t - \mathbb{E}[\overline{X}_t]|^2
+ \frac32|\OV_t|^2 + \frac12\left(\frac{3\lambda_1}m - \frac{\gamma^2}{m^2} \right)|\overline{X}_t - \overline{Y}_t|^2 \\
+ \frac\gamma{2m} \langle \overline{X}_t  - \mathbb{E}[\overline{X}_t], \OV_t\rangle 
+ \frac{\gamma}m \langle \overline{X}_t - \OY_t,\OV_t \rangle  
\end{multline*}
%}
which acts as the Lyapunov functional of the dynamics, 
similar to the strategy introduced in \cite{carrillo2018analytical} for the analysis of CBO methods.

The requirements on the objective function $\E$ are given by conditions \ref{asm:unique} and \ref{asm:growth} from Assumption \ref{asm:E}, and $\E$ is also assumed to be twice differentiable with $\|\nabla ^2\E \|_\infty \leq C_\E$ for some $C_\E>0$. Converge is proved under specific well-posedness assumption on the initial data:

\begin{theorem}[{\cite[Theorem 4]{huang2023pso}}]
Consider the dynamics \eqref{eq:pso_mono_gen} and the above assumptions on $\E$. Moreover, let us assume the well-preparation of the initial datum $\overline{X}_0$ and $\OV_0$ in the sense that $\mu_1,\mu_2>0$ with 
{\small{\begin{align*}
\mu_1&:= \frac{(\lambda_1 + 2\lambda_2)\gamma}{(2m)^2} - \left( \frac{9\lambda_2^2}{\gamma m} + \frac{3\sigma_2^2}{m^2} + \frac{3\lambda_1 \gamma}{4m^2}\right ) 
\frac{12 \exp(-\alpha \min \E)}{\mathbb{E}[\exp(-\alpha \E(\OY_0))]}
\\
\mu_2&:= \frac{(\lambda_1 + \lambda_2)\gamma}{m^2} + \nu\beta \left( \frac{3\lambda_1}{m} + \frac{\gamma^2}{m^2}\right) - \frac{8\nu^2\gamma}{m} - \frac{\lambda_2^2\gamma}{2m^2\lambda_1} - \frac{3\sigma_2^2}{2m^2}\\
&\qquad - \left(\frac{9\lambda_2^2}{\gamma m} + \frac{3\sigma_2^2}{m^2} \right) - \left(\frac{9\lambda_2^2}{\gamma m} + \frac{3\sigma_2^2}{m^2} + \frac{3\lambda_1\gamma}{(2m)^2}\right)\frac{24 \exp(-\alpha \min \E)}{\mathbb{E}[\exp(-\alpha \E(\OY_0))]}
\end{align*}}}
and it holds
{\small{\begin{multline*}
\left(\frac{\alpha \nu m}{\lambda_1 \chi}(C_\E + 2\alpha^2) + \frac{24C^2_\E\nu}{\alpha \chi^3}\right)
\frac{\mathbb{E}[\mathcal{H}(0)]}{\mathbb{E}[\exp(-\alpha (\E(\OY_0)-\min\E))]} \\
+ \frac{6\nu}{\alpha\chi}\frac{\mathbb{E}[|\nabla \E(\overline{X}_0)|^2]}{\mathbb{E}[\exp(-\alpha (\E(\OY_0)-\min\E))]}< \frac{3}{32}
\end{multline*}}}
where 
{\small{\[
\chi:= \frac25 \frac{\min\{\gamma/(2m),\mu_1,\mu_2\}}
{\left((\gamma/(2m))^2 + 1 + 3\lambda_1/m + 2(\gamma/m)^2 \right)}\,.
\]}}
Then $\mathbb{E}[\mathcal{H}(t)]$ converges exponentially fast with rate $\chi$ to $0$ as $t\to\infty$. Moreover, there exists some $\tilde x$, which may depend on $\alpha$ and $f_0$, such that 
\[
\mathbb{E}[\overline{X}_t]\longrightarrow \tilde x, \quad \m^\alpha_\E[\rho^y(t)] \longrightarrow \tilde x \quad \textup{as}\quad t \to \infty
\]
exponentially fast with rate $\chi/2$. 

Eventually, for any given accuracy $\ve>0$, there exists $\alpha_0$, which may depend on the dimension $d$, such that for all $\alpha > \alpha_0$, $\tilde x$ satisfies 
\[\E(\tilde x) - \min \E \leq \ve\,. \]
\end{theorem}

\begin{remark} \label{rmk:inertia}
The inertia mechanism did not belong to the original formulation \cite{kennedy1995,kennedy1997}, but was introduced later in \cite{eberhart1998} as an additional parameter to balance exploration and exploitation in the heuristic strategy. Indeed, large inertia values reduce the contributions given by personal and global bests, favoring the exploration of the search space. This is particularly desirable at the beginning of the computation, and so several adaptive strategies have been proposed for the inertia of the parameter $m$, see \cite{pso2022survey}. %The zero inertia limit of  the kinetic model \eqref{eq:pso:strong} has been derived in \cite{grassi2021pso,grassi2021psosurvey}.
\end{remark}

\subsection{Zero-inertia limit}

As for the kinetic models of SA and GA algorithms, taking suitable scalings of the dynamics allows us to draw a connection between the metaheuristic strategies. In the case of PSO, we consider the generalized model \eqref{eq:pso_mono_gen} in the zero-inertia limit $m\to 0^+$. This will allow us to recover a CBO-type dynamics of first order, with additional memory effects  \cite{borghi2023cbome,Riedl_2024}.
We follow the formal derivation illustrated in \cite{grassi2021pso} and refer to \cite{grassi2021psosurvey} for a rigorous proof.

Let $f = f(x,v,y,t)$ be the law of the SDE system \eqref{eq:pso_mono_gen} which evolves according to 
\begin{multline*}
\partial_t f + v\cdot \nabla_x f + 
+ \nabla_y \cdot \left(\nu(x - y)H^\beta(\E(y) - \E(x))f \right)\\
+ \frac1{m}\nabla_v
\cdot(mv f + \lambda_1(y-x)f +  \lambda_2(\m^\alpha[\rho^y(t)] - x)f)
 = L_m(f)    
\end{multline*}
where we used $\gamma = 1-m$ and introduced 
\[
L_m(f) := \frac1{2m}\sum_{\ell=1}^d 
\Sigma(x_\ell,y_\ell,t)
\partial_{v_\ell}^2
\left( \frac{2fv_\ell}{\Sigma(x_\ell,y_\ell,t)^2} + \frac1{m}\partial_{v_\ell}f\right)
\]
with $\Sigma(x_\ell,y_\ell,t) = \sigma_2^2(m^\alpha[\rho^y(t)]_\ell - x_\ell)^2 + \sigma_1^2(x_\ell - y_\ell)^2$.
Then, we consider the local Maxwellian with unitary mass and zero momentum
\begin{align*}
\mathcal{M}_m(x,v,y,t) &= \prod_{\ell=1}^d M_m(x_\ell, v_\ell,y_\ell,t)\\
M_m(x_\ell,y_\ell,v_\ell,t) &=\frac{m^{1/2}}{\pi^{1/2}|\Sigma(x_\ell,y_\ell,t)|} \exp\left(-\frac{mv_\ell^2}{\Sigma(x_\ell,y_\ell,t)^2} \right)\,.
\end{align*}
This allows to re-write $L_m(f)$ as 
\[L_m(f) = \frac1{2m^2}\sum_{\ell=1}^d 
\Sigma(x_\ell,y_\ell,t)
\partial_{v_\ell}
\left( f \partial_{v_\ell} \log \left(\frac{f}{M_m(x_\ell,y_\ell,v_\ell,t)} \right)\right)\]
and to note that $L_m(f)$ is of order $1/m^2$. Therefore, for $m\ll 1$ we can consider 
\[f(x,v,y,t) = \rho(x,y,t)\mathcal{M}_m(x,y,v,t)\]
and integrate the evolution with respect to $v$, multiply the same equation by $v$ and integrate again with respect to $v$. With such procedure one obtains, see  \cite{grassi2021pso}, the macroscopic system 
\begin{align*}
\partial_t \rho + \nabla_x \cdot(\rho u) + \nabla_y \cdot 
\left(\nu(x-y)H^\beta(\E(y) - \E(x))\rho\right) &= 0 \\
\partial_t (\rho u)_j + \frac{\sigma^2}{2m} \partial_{x_\ell} \left(\Sigma(x_\ell,y_\ell,t)^2 \rho \right) & = \\
-\frac\gamma m (\rho u)_\ell + \frac1m \big( \lambda_1(y_\ell - x_\ell) + &\lambda_1(\m^\alpha[\rho^y(t)]-x_\ell)\big) \rho
\end{align*}
where $\rho u = \int v f\d v$. Finally, by formally taking the zero-inertia limit $m \to 0^+$ we obtain
\begin{multline*}
\partial_t \rho + \nabla_x \cdot \left( (\lambda_1(y-x) + \lambda_2(m^\alpha_\E[\rho^y(t)] - x))\rho\right) + \nabla_y \cdot \left(\nu (x-y)H^\beta(\E(y) - \E(x))\rho \right) \\
 = \frac12\sum_{\ell = 1}^d \partial_{x_\ell}^2 \left( (\sigma_1^2(x_\ell - y_\ell)^2 + \sigma_2^2(m^\alpha_\E[\rho^y(t)]_\ell - x_\ell)^2) \rho\right)\,,
\end{multline*}
which is a CBO dynamics with memory effects and anisotropic noise. Convergence towards global minimizers of this above dynamics has been studied in \cite{borghi2023cbome} and \cite{Riedl_2024} under different assumptions.

\section{Ensemble Kalman Filter}
\label{sec:ekf}
Different from the optimization methods introduced previously, (ensemble Kalman) filtering methods have been applied to solving inverse problems~\cite{evensen1994sequential,herty2019kinetic,herty2020continuous,schillings2017analysis,MR4619819,MR4476008,MR4425851} 
\begin{align}\label{inverse pb}
    y = \mathcal{F}(x) + \eta,
\end{align}
where $\mathcal{F}$ is a possible nonlinear forward operator on $\mathbb{R}^d.$ The unknown state is $x \in \mathbb{R}^d$ and $\eta$ is (observational) noise, usually, assumed to be normally distributed with known co-variance, i.e., $\eta \sim \mathcal{N}(0,\Gamma).$ This problem can be formulated as a weighted least--square problem, being a particular, minimization problem, i.e., 
\begin{align}
    \mathcal{E}(x) = \frac12 | \Gamma^{-\frac12} \left(y - \mathcal{F}(x) \right) |^2.
\end{align}
Contrary to the methods introduced before the structure of $\mathcal{E}(\cdot)$ is relevant in the analysis of the method, see e.g. \cite{aanonsen2009ensemble, blmker2018strongly, blomker2019well, chada2020tikhonov, iglesias2015iterative,janjic2014conservation,schillings2017analysis, schwenzer2020identifying,MR3654885}.    

\subsection{The heuristic strategy}

The ensemble Kalman filter is an iterative filtering method that sequentially updates each member of an ensemble or particle $k=1,\dots, N$ of particles $x_k$ in the state space $\mathbb{R}^d$ by the Kalman update formula.  For $y$ and $\mathcal{F}$ each member $i$ of the ensemble is propagated at the iteration $k$ according to  
\begin{equation}\label{eq:updating0}
    x^{j}_{k+1}=x^{j}_{k} + C(\textbf{x}_k)\left( D(\textbf{x}_k) + \frac{1}{\Delta t} \Gamma^{-1} \right)^{-1}\left[ y - \mathcal{F}(x^j_{k})  \right],
\end{equation}
with initial data $x_0^i$,  where $C(\textbf{x}_k)$ and $D(\textbf{x}_k)$ are the covariance matrices depending on the set  of particles $\textbf{x}_k=(x^1_k, \dots, x_k^N )$ and given by 
  \begin{align}
    \OX_k&:=\frac{1}{N}\sum_{i=1}^N x^j_k, \quad \overline{\mathcal{F}}_k:=\frac{1}{N}\sum_{i=1}^N \mathcal{F}(x^j_k), \notag\\
    C(\textbf{x}_k)&= \frac{1}{N}\sum_{i=1}^N ( x^i_k  -\OX_k)\otimes(\mathcal{F}(x^i_k)-\overline{\mathcal{F}}_k), \label{eq:C} \\ 
    D(\textbf{x}_k)& =\frac{1}{N}\sum_{i=1}^N
    \left[ \mathcal{F}(x^i_k) - \overline{\mathcal{F}}_k \right]\otimes\left[\mathcal{F}(x^i_k) - \overline{\mathcal{F}}_k\right].\label{eq:D}
\end{align}
The algorithm is stated for fixed data $y$, but this could be dependent on $k$. The choice of the initial data $x_0^i$ is arbitrary at the moment, but it has been shown, that (under assumptions) the convex hull of $\{ x^i_k : i=1,\dots,N \}$ (at any time step $k$) is contained in the convex hull of the initial data, see e.g. \cite{MR3654885}. This can be seen as a regularization of the original inverse problem \eqref{inverse pb} through the method. The basic algorithm is given by Algorithm \eqref{alg:enkf}. 

\begin{algorithm}
\caption{Ensemble Kalman Filter}
\label{alg:enkf}
\hspace{-3mm}
\DontPrintSemicolon
\begin{tabular}{ l l l }
$N \in \mathbb{N}$ & number of particles \\  
$y,  \Gamma$ & given data  \\
$\Delta t>0$ & time step
\end{tabular}\tcc*[r]{\scriptsize{parameters}}

$x_{0}^i \in \Rd$ for $i = 1, \dots, N$ \tcc*[r]{\scriptsize{initialization}}
$k = 0$\;

\Repeat{stop condition is satisfied}{
 
Evaluate $\mathcal{F}(x_\k^1),  \dots, \mathcal{F}(x_\k^N)$ \tcc*[r]{\scriptsize{computation of cost}}
$\textbf{x}_k = (x_k^1, \dots, x_k^N)$\;
Compute the co-variance matrices  $C(\textbf{x}_k)$ and $D(\textbf{x}_k)$ with \eqref{eq:C}, \eqref{eq:D}\;
  \For{$i = 1, \dots,  N$}{
$x^i_{k+1} = x^i_k + C(\textbf{x}_k) \left( D(\textbf{x}_k) + \frac{1}{\Delta t} \Gamma^{-1}  \right)^{-1} ( y - \mathcal{F}(x^i_k) ) $  \\
\tcc*[r]{\scriptsize{update particle state}}
}
    $k = k+1$\;
  }
  \KwOut{$\overline{x}_{\k}  = (1/N) \sum_{i=1}^N x^i_{k} $ \tcc*[r]{\scriptsize{return mean value}}}
  \end{algorithm}

\par 
The original formulation has been formulated discrete in iteration (resp. time $k$). However, the limiting equation for $\Delta t \to 0$ under the scaling $\Gamma^{-1}=\Delta t \Gamma^{-1}$ of the previous dynamics has been studied and analyzed, e.g.  \cite{schillings2017analysis,herty2019kinetic}: 
\begin{equation}\label{eq:updating}
    \frac{\d}{\d t} x^{j}_t  =  C(\textbf{x}_t) \Gamma^{-1} \left[ y - \mathcal{F}(x^j_t)  \right]. 
\end{equation}
In order to gain a heuristic understanding, the linear case has been studied, i.e., $\mathcal{F}(x)= {\bf F}x$ for a given matrix ${\bf F}\in \mathbb{R}^{d\times d}$. In this case, a direct computation shows that the  dynamics \eqref{eq:updating} reads
\begin{equation}
\begin{split}
    \frac{\d}{\d t} x^{j}_t &=  - \overline{C}(\textbf{x}_t) \nabla_x \mathcal{E}(x^j), \; x^{j}(0)=x_{j,0}, \\
    \overline{C}(\textbf{x}_t) &= \frac{1}{N}\sum_{i=1}^N ( x^i_t  -\overline{x}_t) \otimes
    ( x^i_t  -\overline{x}_t). 
    \end{split}\label{eq:updating2}
\end{equation}
Since $\overline{C}(\textbf{x}_t)$ is positive definite, the dynamics of the ensemble Kalman filter can be understood as a (preconditioned) gradient descent method. This is seen as a heuristic argument for why Kalman filtering performs well on a variety of inverse problems of the type \eqref{inverse pb}, see e.g., \cite{MR3764752} for theoretical and numerical analysis in the linear, but noisy case.
Note that the formulation \eqref{eq:updating2} is misleading since the original method \eqref{eq:updating} is derivative free. The Ensemble Kalman method is therefore a method where  $\mathcal{E}$ is non-differentiable similar to particle swarm methods or genetic algorithms. 
\par 
The dynamics \eqref{eq:updating2} has been analysed with respect to convergence properties. The right--hand side of equation \eqref{eq:updating2} is Lipschitz continuous yielding the existence of a unique solution $x_t^j$ for $j=1,\dots,N.$ Define the ensemble spread $e^j_t$ 
\begin{align}
    e^j_t = x^j_t - \bar{x}_t, \; r^j_t = x^j_t - x^*.
\end{align}
Then the following result \cite{schillings2017analysis} gives sufficient conditions for a collapse towards
the mean value at an algebraic rate. 

\begin{theorem}\cite[Theorem 3]{schillings2017analysis} 
Assume that $y={\bf F}x^*$ for some $x^* \in \R^d.$ Let $\Gamma$ be positive definite and define $E_t \in \R^{N\times N}$  by 
\begin{align}
(E_t)_{k,\ell} = e^\ell_t {\bf F}^T \Gamma {\bf F} e^k_t.     
\end{align}
Then, $\| E_ t \| \to 0$ at rate $\mathcal{O}(N t^{-1})$ for $t \to \infty.$
\end{theorem}

The convergence of the residual $r^j_t$ has been shown~\cite[Proposition 3.4]{armburster} for a slightly modified version of the dynamics compared with  \eqref{eq:updating2} and for parameters $\kappa,\beta \in \R$ given by 
\begin{equation}
\begin{split}
    \frac{\d}{\d t} x^{j}_t &=  - \overline{C}(\textbf{x}_t;\kappa) \left( \nabla_x  \mathcal{E}(x^j) + \beta (\bar{x}_t - x^j_t \right), \; x^{j}(0)=x_{j,0}, \\
    \overline{C}(\textbf{x}_t) &= \frac{1}{N}\sum_{i=1}^N ( x^i_t  - \kappa \overline{x}_t) \otimes
    ( x^i_t  - \kappa \overline{x}_t). 
    \end{split}\label{eq:updating3}
\end{equation}

\begin{theorem}\cite[Proposition 3.4]{armburster}
Assume that $y={\bf F}x^*$ for some $x^* \in \R^d.$ Assume $x^*$ is a Karush--Kuhn--Tucker point of the minimization problem $\min_u \mathcal{E}(u).$ 
Then, for $\kappa(2-\kappa)<1$ and $\beta<0$ we have 
\begin{align}
    \lim\limits_{t\to\infty} \| r^j_t \|^2 = 0. 
\end{align}
\end{theorem}
The modified dynamics \eqref{eq:updating3} is motivated by a study of the moments \eqref{eq:moments} of the mean-field model \eqref{eq:meanfieldENKF}
below. 

\subsection{The mean-field model}

For the time continuous formulation \eqref{eq:updating} and $N\to \infty$ a mean-field limit is obtained. Rigorous results can be found e.g. in \cite{herty2019kinetic} and in \cite{diez2022review1,diez2022review2} for general mean-field results on interacting particle systems. From a kinetic perspective, the particle dynamics \eqref{eq:updating} is simpler compared with particle swarming dynamics. This is due to the fact, that the interaction of particles is only due to the operator $C(\textbf{x}_t).$
The limit for $N\to\infty$ can therefore be understood as a mean-field limit of the particle dynamics. Convergence of the empirical measure $f^N(t) = (1/N) \sum_{i=1}^N \delta_{x_i(t)}$ to the solution of the mean-field equation \eqref{eq:meanfieldENKF} $f(t) \in \mathcal{P}(\R^d)$ in 2-Wasserstein distance has been shown, see e.g. \cite{herty2019kinetic}. The propagation of chaos, similar to \eqref{eq:psochaos} or \eqref{eq:chaosga} is also assumed here. The convergence is understood  in the sense introduced in equation \eqref{eq:convergenceW}.
\par 
Since the derivation is similar to the results presented in the section \eqref{sec:pso:derivation}, we give only
the resulting mean-field equation in strong form:
\begin{align}\label{eq:meanfieldENKF}
    \partial_t f + \nabla_x \cdot \left( \mathcal{C}[f(t)] \Gamma^{-1} \left( y - \mathcal{F}(x) \right) f \right) = 0,
\end{align}
with initial data $f(0) = f_0$ for some $f_0\in\mathcal{P}(\Rd)$.
The initial data $f_0$ is obtained as limit in the sense $W_2(f^N_0, f_0) \to 0$ for $N\to \infty.$ 
The operator $\mathcal{C}[\cdot]$  is nonlocal and it can be written in terms of moments. For any $g\in \mathcal{P}(\Rd)$, it is defined as 
\begin{equation}\label{eq:moments}
\begin{split}
    \m[g] := \int x  g(\d x), \quad \m_\mathcal{F}[g] := \int \mathcal{F}(x) g(\d x), \\
    \mathcal{C}[g] := \int  ( x - \m[g] ) \otimes ( \mathcal{F}(x) - 
    \m_\mathcal{F}[g]  ) g(\d x).
\end{split} 
\end{equation}

In the linear case, i.e., $\mathcal{F}(x)={\bf F}(x)$, the nonlocal operator $\mathcal{C}[g]$ further simplifies. The operator can be written explicitly in terms of $(\m, \textup{E} )$ where $\textup{E}$ is defined as second moment of $g,$ i.e., 
\begin{align*}
    \textup{E}[g] = \int x \otimes x \; g(\d x).
\end{align*}
The time evolution of $(\m,\textup{E})$ can be written as a closed system of ordinary differential equations, but it has been shown that even in the case $d=1$, the stationary states $(\m,\m^2)$ for $\m\in\R$ are  non–hyperbolic Bogdanov–Takens–type equilibria \cite{armburster}. Hence, they are  not asymptotically stable. This has been the motivation to consider modified dynamics as e.g. \eqref{eq:updating3}. In this case, the stationary states $(\m,\m^2)$ of the moment system corresponding to the  mean-field limit of the modified dynamics \eqref{eq:updating3} are globally asymptotically stable for $\kappa(2-\kappa)<1$ and $\beta<0$, see \cite[Proposition 3.2]{armburster}.

\begin{remark}
Only the basic method presented by equation \eqref{eq:updating} has been presented. The convergence of the numerical method has also been studied \cite{MR3840898}. Many extensions towards inverse problems with additional box constraints \cite{MR4617850}, or linear constraints \cite{MR4121318}, multi-scale \cite{MR4184278}, or multi-objective optimization \cite{MR4605205} have been proposed. The interaction kernel $C(\textup{\textbf{x}}_t)$ in the update \eqref{eq:updating}, if done naively, requires a quadratic in $d$, which may make the method computationally infeasible. However, it has been shown to be efficient even for the training of neural networks, see e.g.,  \cite{MR4409716,alper}. Also, efficient sub-sampling strategies have been proposed to reduce the computational overhead, see e.g., \cite{MR4629223}.
\end{remark}

%\newpage 
%\section{Conclusions}
%{\color{blue} TBC }

\section{Where we are and what's next}

This review represents an initial attempt to present some of the main metaheuristic optimization algorithms—Simulated Annealing, Genetic Algorithms, Particle Swarm Optimization, and the Ensemble Kalman Filter—in a unified framework based on the principles of kinetic theory, providing a probabilistic description of their associated dynamics. These algorithms are essential tools in scenarios where traditional optimization methods fall short, particularly when dealing with non-convex, non-differentiable, or black-box objective functions.

By employing this formalism, we establish a rigorous mathematical foundation for understanding the dynamics of these methods, revealing their underlying mechanisms within a coherent framework. The models derived not only facilitate the analysis of convergence properties using classical PDE techniques but also enable the development of new optimization strategies informed by the insights gained from kinetic descriptions. This connection to kinetic theory fosters a deeper understanding of the relationships among various metaheuristic methods, providing a pathway for future innovations in optimization techniques.

%Reflecting on our current understanding, it is evident that integrating kinetic theory into the analysis of metaheuristic optimization has paved the way for further research. The promise of adaptive methods informed by kinetic models holds significant potential for enhancing the efficiency and effectiveness of optimization strategies across a range of applications.

Looking ahead, several promising directions for future research emerge. Continued exploration of kinetic models will likely yield deeper insights into the behavior of metaheuristic algorithms, enabling the development of adaptive methods that can better navigate complex fitness landscapes or operate effectively in uncertain environments. Moreover, the ability to establish convergence guarantees under mild assumptions encourages further investigation into the theoretical foundations of these algorithms, potentially leading to the design of new heuristics grounded in rigorous mathematical principles.

The framework established in this review also paves the way for the analysis of additional metaheuristics, enabling the discovery of further analogies and potential improvements. By bridging the gap between heuristic optimization and rigorous mathematical foundations, we not only deepen our understanding of existing methods but also set the stage for future innovations capable of addressing increasingly complex optimization challenges

%------
% Insert acknowledgments and information
% regarding funding at the end of the last
% section, i.e., right before the bibliography.
%------

%begin{ack}
%We thank X.
%\end{ack}

\begin{funding}
This work has been written within the activities of GNCS group of INdAM (Italian National Institute of High Mathematics). The research of G.B. and L.P. has been supported by the Royal Society under the Wolfson Fellowship “Uncertainty quantification, data-driven simulations and learning of multiscale complex systems governed by PDEs”. L.P. also acknowledges the  partial support 
by European Union -
NextGenerationEU through the Italian Ministry of University and Research as
part of the PNRR – Mission 4 Component 2, Investment 1.3 (MUR Directorial
Decree no. 341 of 03/15/2022), FAIR “Future” Partnership Artificial Intelligence Research”, Proposal Code PE00000013 - CUP DJ33C22002830006) and by MIUR-PRIN Project 2022, No. 2022KKJP4X “Advanced numerical methods for time dependent parametric partial differential equations with applications”.
M.H. thanks the Deutsche Forschungsgemeinschaft (DFG, German Research Foundation) for the financial support through 320021702 / GRK2326, and for the financial support through 442047500/SFB1481 within the projects B04 (Sparsity f{\"o}rdernde Muster in kinetischen Hierarchien), B05 (Sparsifizierung zeitabh{"a}ngiger Netzwerkflußprobleme mittels diskreter Optimierung) and B06 (Kinetische Theorie trifft algebraische Systemtheorie). M.H. and L.P. received funding from the European Union’s Horizon Europe research and innovation programme under the Marie Sklodowska-Curie Doctoral Network Datahyking (Grant No. 101072546).

\end{funding}

%------
% Insert the bibliography.
%------
% \bibliographystyle{ems}
  \bibliographystyle{abbrv}
 \bibliography{bibfile,references}
%\begin{thebibliography}{99}

%------ Example for a paper in journal:
% \bibitem{article1}
% A.~Petrunin, Parallel transportation for Alexandrov space with curvature bounded below.
% \emph{Geom. Funct. Anal.} \textbf{8} (1998), no.~1, 123--148.
% \Zbl{0903.53045} \MR{1601854}

%------ Example for a book:
% \bibitem{book1}
% W.~P. Ziemer, \emph{Weakly differentiable functions}.
% Grad. Texts in Math. 120,  Springer, New York, 1989.
%\Zbl{0692.46022} \MR{1014685}

%------ Example for a paper in a book:
% \bibitem{incollection1}
% J.~S. Milne, Introduction to Shimura varieties.
% In \emph{Harmonic analysis, the trace formula, and Shimura varieties},
% edited by M.~W. Marcellin and E.~Giorgi, pp. 265--378,
% Clay Math. Proc. 4, Amer. Math. Soc., Providence, RI, 2005.
% \Zbl{1148.14011} \MR{2192012}

%------ Example for a preprint on arXiv:
% \bibitem{preprint1}
% D.~V. Nguyen, S.~K. Chilappagari, M.~W. Marcellin, and B.~Vasic,
% LDPC codes from latin squares free of small trapping sets,
% 2010, \href{http://arxiv.org/abs/1008.4177}{arXiv:1008.4177}.

%------ Example for a report:
% \bibitem{report1}
% J.~Schöberl, Commuting quasi-interpolation operators.
% Technical report isc-01-10-math, Texas A\&M University, 2001,
% \url{www.isc.tamu.edu/publications-reports/tr/0110.pdf}.

%------ Example for a thesis:
% \bibitem{thesis1}
% E.~Giorgi, \emph{The geometric universe}.
% Ph.D. thesis, University of Maryland, College Park, 2002.

%\end{thebibliography}

\end{document}